\documentclass{article}

\newif\ifarxiv
\arxivtrue

\usepackage{graphicx}
\usepackage{amsmath,amssymb,amsthm,bm}
\usepackage[a4paper]{geometry}
\usepackage{algorithm,algorithmicx,algpseudocode}
\newcommand*{\SComment}[1]{\Comment{{\scriptsize #1}}} 
\usepackage[colorlinks=false,pdfborder={0 0 0}]{hyperref}
\usepackage{authblk}
\usepackage{color}
\usepackage{comment}
\usepackage{xspace}
\usepackage{subfigure}

\graphicspath{{fig/}}
\DeclareGraphicsExtensions{.mps,.pdf,.eps,.png}

\newenvironment{keywords}{\medskip\textbf{Keywords:}}{}

\newenvironment{MSCcodes}{\medskip\textbf{MSC codes:}}{}

\newtheorem{theorem}{Theorem}

\newtheorem{lemma}{Lemma}

\newtheorem{corollary}{Corollary}
\theoremstyle{plain}

\DeclareMathOperator{\Span}{span}

\newcommand{\fro}{\mathsf F}

\newcommand*{\trans}{^{\top}}

\newcommand*{\itrans}{^{-\top}}

\newcommand{\inv}{{-1}}

\newcommand{\bmat}[1]{\begin{bmatrix}#1\end{bmatrix}}

\newcommand*{\norm}[1]{\Vert#1\rVert}
\newcommand*{\normF}[1]{\bigl\Vert#1\bigr\rVert_{\fro}}

\def\adots{\mathinner{\mkern2mu\raise1pt\hbox{.}\mkern2mu
    \raise4pt\hbox{.}\mkern2mu\raise7pt\hbox{.}\mkern1mu}}

\newcommand{\sigmin}{\sigma_{\min}}
\newcommand{\sigmax}{\sigma_{\max}}

\newcommand{\vD}{\bm{D}}

\newcommand{\vG}{\bm{G}}

\newcommand{\vJhat}{\hat{\bm{J}}}

\newcommand{\vQ}{\bm{Q}}
\newcommand{\vQhat}{\widehat{\vQ}}

\newcommand{\vS}{\bm{S}}

\newcommand{\vU}{\bm{U}}

\newcommand{\vUhat}{\hat{\vU}}
\newcommand{\vV}{\bm{V}}

\newcommand{\vVhat}{\widehat{\bm{V}}}
\newcommand{\vW}{\bm{W}}

\newcommand{\vWhat}{\widehat{\bm{W}}}
\newcommand{\vX}{\bm{X}}

\newcommand{\vY}{\bm{Y}}

\newcommand{\vZ}{\bm{Z}}
\newcommand{\vZhat}{\widehat{\vZ}}

\newcommand{\Phat}{\widehat{P}}

\newcommand{\Rhat}{\widehat{R}}

\newcommand{\Shat}{\hat{S}}

\newcommand{\That}{\widehat{T}}

\newcommand{\Omhat}{\hat{\Omega}}

\newcommand{\Yhat}{\widehat{Y}}

\newcommand{\RR}{\mathcal{R}}
\newcommand{\RRhat}{\widehat{\RR}}
\renewcommand{\SS}{\mathcal{S}}
\newcommand{\SShat}{\widehat{\SS}}

\newcommand{\YY}{\mathcal{Y}}
\newcommand{\YYhat}{\widehat{\mathcal{Y}}}

\newcommand{\bDD}{\bm{\mathcal{D}}}

\newcommand{\bQQ}{\bm{\mathcal{Q}}}
\newcommand{\bQQhat}{\widehat{\bQQ}}
\newcommand{\bQQtil}{\widetilde{\bQQ}}

\newcommand{\bXX}{\bm{\mathcal{X}}}

\newcommand{\DeltaC}{\Delta C}
\newcommand{\DeltaCone}{\DeltaC^{(1)}}
\newcommand{\DeltaCtwo}{\DeltaC^{(2)}}

\newcommand{\DeltaF}{\Delta F}

\newcommand{\DeltaFone}{\DeltaF^{(1)}}
\newcommand{\DeltaFtwo}{\DeltaF^{(2)}}

\newcommand{\DeltaJ}{\Delta J}

\newcommand{\DeltaP}{\Delta P}

\newcommand{\DeltaT}{\Delta T}

\newcommand{\DeltaOmega}{\Delta\Omega}

\newcommand{\DeltavD}{\Delta\vD}

\newcommand{\DeltavG}{\Delta\vG}
\newcommand{\DeltavGone}{\DeltavG^{(1)}}
\newcommand{\DeltavGtwo}{\DeltavG^{(2)}}

\newcommand{\DeltavS}{\Delta\vS}

\newcommand{\DeltavV}{\Delta\vV}

\newcommand{\DeltavVhat}{\Delta\vVhat}
\newcommand{\DeltavW}{\Delta\vW}

\newcommand{\DeltavY}{\Delta\vY}

\newcommand{\DeltabDD}{\Delta\bDD}

\newcommand{\deltaJ}{\delta_{J}}

\newcommand{\deltacholone}{\delta_{\chol_1\xspace}}

\newcommand{\deltacholtwo}{\delta_{\chol_2\xspace}}
\newcommand{\deltaP}{\delta_{P}}
\newcommand{\deltaQ}{\delta_{Q}}
\newcommand{\deltaQS}{\delta_{QS}}

\newcommand{\deltaQU}{\delta_{Q\trans U}}
\newcommand{\deltaQX}{\delta_{Q\trans X}}
\newcommand{\deltaQY}{\delta_{QY}}
\newcommand{\deltaS}{\delta_{S}}
\newcommand{\deltaSS}{\delta_{S\trans S}}

\newcommand{\deltaU}{\delta_{U}}

\newcommand{\deltaUU}{\delta_{U\trans U}}
\newcommand{\deltaX}{\delta_{X}}
\newcommand{\deltaXX}{\delta_{X\trans X}}

\newcommand{\deltaYY}{\delta_{Y\trans Y}}
\newcommand{\deltaZ}{\delta_Z}

\newcommand{\BCGSIRO}{\texttt{BCGSI+}\xspace}	
\newcommand{\BCGSIROA}{\texttt{BCGSI+A}\xspace}	
\newcommand{\BCGSIROAPIPIROone}{\texttt{BCGSI+P-1S-2S}\xspace}	
\newcommand{\BCGSIROAone}{\texttt{BCGSI+A-1s}\xspace}	
\newcommand{\BCGSPIPIRO}{\hyperref[alg:BCGSPIPIRO]{\texttt{BCGS-PIPI+}}\xspace}	
\newcommand{\BCGSPIPIROone}{\hyperref[alg:BCGSPIPIRO1S]{\texttt{BCGSI+P-1S}}\xspace}
\newcommand{\BCGSPIPIROtwo}{\hyperref[alg:BCGSPIPIRO2S]{\texttt{BCGSI+P-2S}}\xspace}

\newcommand{\IOA}[1]{\texttt{IO}_{\mathrm{A}}\left(#1\right)}
\newcommand{\IOAnoarg}{\texttt{IO}_{\mathrm{A}}\xspace}
\newcommand{\IOone}[1]{\texttt{IO}_1\left(#1\right)}
\newcommand{\IOonenoarg}{\texttt{IO}_1\xspace}
\newcommand{\IOtwo}[1]{\texttt{IO}_2\left(#1\right)}
\newcommand{\IOtwonoarg}{\texttt{IO}_2\xspace}

\newcommand{\MGS}{\texttt{MGS}\xspace}	
\newcommand{\HouseQR}{\texttt{HouseQR}\xspace}	
\newcommand{\TSQR}{\texttt{TSQR}\xspace}	
\newcommand{\CholQR}{\texttt{CholQR}\xspace}	

\newcommand*{\chol}{\texttt{chol}\xspace}

\newcommand*{\macheps}{\bm u}

\newcommand*{\bigO}{O}

\newcommand{\monomial}{\texttt{monomial}\xspace}

\newcommand{\default}{\texttt{default}\xspace}
\newcommand{\glued}{\texttt{glued}\xspace}
\newcommand{\piled}{\texttt{piled}\xspace}

\definecolor{plotred}{RGB}{216.7500, 82.8750, 24.9900}
\newcommand{\diff}[1]{\textcolor{black}{#1}}

\newcommand*{\krylov}{\mathcal{K}}
\newcommand*{\kn}{n^*}
\newcommand*{\ki}{k^*}
\newcommand*{\kj}{j^*}

\title{A stable one-synchronization variant of reorthogonalized block classical Gram--Schmidt}
\author[1]{Erin Carson}
\author[1]{Yuxin Ma}
\affil[1]{Department of Numerical Mathematics, Faculty of Mathematics and Physics, Charles University, Sokolovsk\'{a} 49/83, 186 75 Praha 8, Czechia

{\tt Email: carson@karlin.mff.cuni.cz, yuxin.ma@matfyz.cuni.cz}}

\begin{document}
\maketitle

\begin{abstract}
The block classical Gram--Schmidt (BCGS) algorithm and its reorthogonalized variant are widely-used methods for computing the thin QR factorization of a given matrix \(\bXX\) with block vectors due to their lower communication cost compared to other approaches such as modified Gram--Schmidt and Householder QR.
To further reduce communication, i.e., synchronization, there has been a long ongoing search for a variant of reorthogonalized BCGS variant that achieves \(\bigO(\macheps)\) loss of orthogonality while requiring only \emph{one} synchronization point per block column, where \(\macheps\) represents the unit roundoff. 
Utilizing Pythagorean inner products and delayed normalization techniques, we propose the first provably stable one-synchronization reorthogonalized BCGS variant, demonstrating that it has \(\bigO(\macheps)\)  loss of orthogonality under the condition \(\bigO(\macheps) \kappa^2(\bXX) \leq 1/2\), where $\kappa(\cdot)$ represents the condition number. 

By incorporating one additional synchronization point, we develop a two-synchronization reorthogonalized BCGS variant which maintains \(\bigO(\macheps)\) loss of orthogonality under the improved condition \(\bigO(\macheps) \kappa(\bXX) \leq 1/2\).
An adaptive strategy is then proposed to combine these two variants, ensuring \(\bigO(\macheps)\) loss of orthogonality while using as few synchronization points as possible under the less restrictive condition \(\bigO(\macheps) \kappa(\bXX) \leq 1/2\).
As an example of where this adaptive approach is beneficial, we show that using the adaptive orthogonalization variant, \(s\)-step GMRES achieves a backward error comparable to \(s\)-step GMRES with \BCGSIRO, also known as \texttt{BCGS2}, both theoretically and numerically, but requires fewer synchronization points.


\end{abstract}

\section{Introduction}
\label{sec:introduction}
Given a matrix \(\bXX \in \mathbb{R}^{m\times n}\), this work considers computing the thin QR factorization for \(\bXX\), that is,
\[
    \bXX = \bQQ \RR,
\]
where \(\bQQ \in \mathbb{R}^{m\times n}\) is orthonormal and \(\RR\in \mathbb{R}^{n \times n}\) is upper triangular.
The classical Gram--Schmidt algorithm is a widely-used method for this problem.

Reducing communication, i.e., data movement and synchronization, has become increasingly crucial for the performance of traditional solvers. Thus blocking, which reduces synchronization points and leverages BLAS-3 operations, is an appealing approach; see the survey~\cite{CLRT2022}.
Consequently, the block classical Gram--Schmidt (BCGS) algorithm and its variants have attracted significant interest.
To simplify, we consider \(n = ps\) and let \(\bXX \in \mathbb{R}^{m\times ps}\) consist of \(p\) block vectors, such that \(\bXX = \bmat{\vX_1& \dotsi & \vX_p}\) where \(\vX_k \in \mathbb{R}^{m\times s}\) for any \(k = 1, 2, \dotsc, p\).
BCGS is employed, for example, in communication-avoiding Krylov methods~\cite{Carson2015,Hoemmen2010,YTHBSE2020}, which have been demonstrated to outperform non-block methods in many practical scenarios~\cite{Hoemmen2010,YATHD2014}.
In BCGS, particularly when integrated within communication-avoiding Krylov methods, we typically focus on the backward error
\[
\frac{\norm{\bXX - \bQQhat \RRhat}}{\norm{\bXX}},
\]
and the loss of orthogonality (LOO)
\[
\norm{I - \bQQhat\trans \bQQhat},
\]
where \(\hat{\cdot}\) denotes computed quantities and \(\norm{\cdot}\) denotes the \(2\)-norm. 
These two aspects will directly impact the stability of communication-avoiding Krylov subspace methods and other downstream applications.
It should be noted that the backward error is \(\bigO(\macheps)\) for all commonly used BCGS methods, where \(\macheps\) represents the unit roundoff, making the LOO typically the primary concern.
Moreover, in the paper, we use \(\kappa(\bXX) = \sigmax(\bXX)/\sigmin(\bXX)\) to represent the 2-norm condition number if \(\sigmin(\bXX)>0\), where \(\sigmax(\bXX)\) and \(\sigmin(\bXX)\) denote the largest and smallest singular value of \(\bXX\), respectively.

Considering a distributed memory parallel setting, we define a synchronization point as a global communication involving all processors, such as \texttt{MPI\_Allreduce}.
The traditional BCGS method involves two synchronization points (i.e., block inner products and \emph{intraorthogonalization} which orthogonalize vectors within a block column) per iteration (i.e., per block column) but suffers from instability, as its LOO is not \(\bigO(\macheps)\)~\cite{CLMO2024-ls}.
To improve its LOO, one can use a reorthogonalization technique, essentially running the BCGS iteration twice in the for-loop. This reorthogonalized variant  \BCGSIRO, also known as \texttt{BCGS2}, achieves \(\bigO(\macheps)\) LOO under certain conditions regarding the condition number \(\kappa(\bXX)\), as analyzed in~\cite{Barlow2024,BS2013,CLMO2024-ls}.
However, \BCGSIRO's main drawback is that it requires four synchronization points per iteration.

Due to the expensive nature of synchronization, there has been significant interest in minimizing the number of synchronization points per iteration.
Synchronization is typically needed for computing inner products or norms of large matrices.
One possible strategy to reduce the number of synchronization points is to organize these inner product and norm calculations as collectively as possible.
As proposed in~\cite{CLMO2024-P}, \BCGSPIPIRO is a variant of reorthogonalized BCGS that utilizes Pythagorean inner products to reduce synchronization points from four to two, achieving \(\bigO(\macheps)\) LOO under the assumption \(\bigO(\macheps) \kappa^2(\bXX) \leq 1\).
Derived in~\cite{CLMO2024-ls}, \BCGSIROAone, which resembles \texttt{BCGS+LS} proposed in~\cite{YTHBSE2020}, is a one-sync reorthogonalized BCGS variant created by removing one intraorthogonalization and delaying another.
Additionally, it is the block version of \texttt{DCGS2}~\cite{BLTSYB2022} or \texttt{CGS-2}~\cite{SLAYT2021}, which is demonstrated to have \(\bigO(\macheps)\) LOO in~\cite{CLMO2024-ls}.
Unfortunately, \BCGSIROAone can only achieve \(\bigO(\macheps) \kappa^2(\bXX)\) LOO given the assumption \(\bigO(\macheps) \kappa^3(\bXX) \leq 1\)~\cite{CLMO2024-ls}.

In this paper, our aim is first to formulate a stable one-sync reorthogonalized BCGS in which the backward error and LOO remain at the level \(\bigO(\macheps)\).
In~\cite{CLMO2024-ls}, a two-sync variant reorthogonalized BCGS, called \texttt{BCGSI+A-2s}, stems from removing the first intraorthogonalization and using Pythagorean-based Cholesky QR as the second intraorthogonalization.
Unlike \BCGSPIPIRO, which also features two synchronization points per iteration, this version only reaches \(\bigO(\macheps) \kappa^2(\bXX)\) LOO.
According to the analysis in~\cite{CLMO2024-ls}, comparing these two variants shows that the first intraorthogonalization plays a significant role in the LOO.

Consequently, we reintegrate the omitted intraorthogonalization into \BCGSIROAone, which can be seen either as utilizing Pythagorean inner products to boost the LOO of \BCGSIROAone or applying the delayed intraorthogonalization concept to \BCGSPIPIRO, to eliminate an additional synchronization point.
Both approaches lead to the one-sync reorthogonalized BCGS method, denoted as \BCGSPIPIROone in this paper, which achieves \(\bigO(\macheps)\) LOO, provided \(\bigO(\macheps) \kappa^2(\bXX) \leq 1\). 

Moreover, observing that the assumption \(\bigO(\macheps) \kappa^2(\bXX) \leq 1\) in LOO arises from employing Pythagorean-based Cholesky QR for the first intraorthogonalization, we allow for the possibility of a more stable intraorthogonalization method, which requires at least one extra synchronization.
Tall-Skinny QR (TSQR) is a common choice for intraorthogonalization in communication-avoiding methods like low-synchronization BCGS and \(s\)-step GMRES, since it requires only one synchronization~\cite{Hoemmen2011}.
This results in \BCGSPIPIROtwo with \(\bigO(\macheps)\) LOO, provided \(\bigO(\macheps) \kappa(\bXX) \leq 1\).
This new variant of BCGS maintains the same LOO property as \BCGSIRO while reducing the synchronization needs from four to two points per iteration.

For the one-sync variant \BCGSPIPIROone, the requirement on the condition number of the input matrix $\bXX$ makes this approach generally unsuitable for use within communication-avoiding Krylov subspace methods.
For example, communication-avoiding GMRES, also known as \(s\)-step GMRES, analyzed in~\cite{CM2024}, requires an orthogonalization method yielding a well-conditioned \(Q\)-factor under the assumption \(\bigO(\macheps) \kappa^{\alpha}(\bXX) \leq 1\), with \(\alpha = 0\) or \(1\) in order to guarantee a backward error at the \(\bigO(\macheps)\) level.
Thus, we suggest a new adaptive approach to orthogonalization within communication-avoiding GMRES: primarily adopting \BCGSPIPIROone and switching to \BCGSPIPIROtwo if \(\bigO(\macheps) \kappa^2(\bXX_k) > 1\), where \(\bXX_k = \bmat{\vX_1& \vX_2& \dotsi & \vX_k}\).
This switching condition, i.e., \(\bigO(\macheps) \kappa^2(\bXX_k) > 1\), can be checked using \(\bigO(ns^2)\) operations without necessitating extra synchronization, according to our analysis.

The remainder of this paper is organized as follows.
We propose the one-sync reorthogonalized BCGS, and analyze its backward error and LOO in Section~\ref{sec:algo-onesync}.
Then in Section~\ref{sec:algo-twosync}, we propose a two-sync reorthogonalized BCGS and show that its LOO property is as same as \BCGSIRO.
In Section~\ref{sec:orthoGMRES}, we demonstrate how to employ these low-sync reorthogonalized BCGS in \(s\)-step GMRES.
Numerical experiments are presented in Section~\ref{sec:experiments} to compare the new variants introduced in Sections~\ref{sec:algo-onesync}, \ref{sec:algo-twosync}, and~\ref{sec:orthoGMRES} with \BCGSIROA, \BCGSIROAone, and \BCGSPIPIRO, and to further compare these variants employed in \(s\)-step GMRES.

We introduce some notation used throughout the paper before proceeding.
Uppercase bold Roman scripts, for example, \(\bXX\) and \(\bQQ\), denote the matrices containing \(p\) block vectors.
Each block vector is denoted by uppercase Roman letters \(\vX_k\) and \(\vQ_k\), i.e., \(\bXX_k = \bmat{\vX_1& \vX_2& \dotsi & \vX_k}\) and \(\bQQ_k = \bmat{\vQ_1& \vQ_2& \dotsi & \vQ_k}\).
The regular scripts (without using bold and Roman styles) represent for the column of matrices, i.e., \(X_i\) denotes the \(i\)-th column of the matrix \(\bXX\).
For square matrices, we use uppercase Roman scripts, for example, \(\RR\), \(\SS\), and \(\YY\), to represent \(ps \times ps\) matrices.
We use MATLAB indexing to denote submatrices.
For example,
\[
    \RR_{1:k-1,k} = \bmat{R_{1,k}\\ R_{2,k}\\ \vdots \\ R_{k-1,k}}, \quad
    \SS_{1:k-1,k} = \bmat{S_{1,k}\\ S_{2,k}\\ \vdots \\ S_{k-1,k}}, \quad\text{and}\quad
    \YY_{1:k-1,k} = \bmat{Y_{1,k}\\ Y_{2,k}\\ \vdots \\ Y_{k-1,k}}.
\]
In addition, \(\RR_{kk} = \RR_{1:k, 1:k}\), \(\SS_{kk} = \SS_{1:k, 1:k}\), and \(\YY_{kk} = \YY_{1:k, 1:k}\).
We also use \(\normF{\cdot}\) to denote the Frobenius norm.
The functions $[\vQ, R] = \IOA{\vX}$, $[\vQ, R] = \IOone{\vX}$, and $[\vQ, R] = \IOtwo{\vX}$ represent intraorthogonalization, which are methods to perform the thin QR factorization of a block vector $\vX$.
Algorithms such as Householder QR (\HouseQR), TSQR (\TSQR), classical Gram--Schmidt, modified Gram--Schmidt (\MGS), or Cholesky QR (\CholQR), described in~\cite{CLRT2022}, can be applied as intraorthogonalization routines.
\section{One-sync Pythagorean-based BCGS}
\label{sec:algo-onesync}
In this section, we first introduce a new one-sync reorthogonalized BCGS method, and then analyze the LOO of this method. Our analysis will show that the LOO is bounded by \(\bigO(\macheps)\) LOO under the condition \(\bigO(\macheps) \kappa^2(\bXX) \leq 1/2\).

\subsection{Algorithm}
In \cite{CLMO2024-P}, the authors analyze a reorthogonalized variant of BCGS using the Pythagorean inner product, referred to as \(\BCGSPIPIRO\) shown in Algorithm~\ref{alg:BCGSPIPIRO}.
They demonstrate that it is guaranteed to achieve \(\bigO(\macheps)\) loss of orthogonality, provided that \(\bigO(\macheps) \kappa^2(\bXX) \leq 1/2\).
It should be noted that there are two synchronization points in Algorithm~\ref{alg:BCGSPIPIRO}:
one is required in Line~\ref{line:bcggspipiro:SS} and the other is required in Line~\ref{line:bcggspipiro:YY-Omega}.

Based on the insights from~\cite{YTHBSE2020}, one strategy to eliminate a synchronization point involves reorganizing the inner products related to large matrices as compactly as possible.
One potential approach is to move Line~\ref{line:bcggspipiro:YY-Omega} immediately after Line~\ref{line:bcggspipiro:SS}.
However, since in Line~\ref{line:bcggspipiro:YY-Omega}, \(\vU_{k+1}\) relies on the results computed in Line~\ref{line:bcgspipiro:Sdiag}, this reordering is not feasible.
Alternatively, Line~\ref{line:bcggspipiro:SS} in the \(k\)-th iteration can be shifted just after Line~\ref{line:bcggspipiro:YY-Omega} in the (\(k-1\))-st iteration. 
Observing the \(k\)-th iteration, it is clear that in Line~\ref{line:bcggspipiro:SS}, \(\vX_{k+1}\trans \vX_{k+1}\) is independent of the computations from Lines~\ref{line:bcggspipiro:Ydiag}--\ref{line:bcgspipiro:Rdiag} in the (\(k-1\))-st iteration, and
\begin{equation*}
    \begin{split}
        \bQQ_{k}\trans \vX_{k+1}
        &= \bmat{\bQQ_{k-1}\trans \vX_{k+1}& \vQ_k\trans \vX_{k+1}},
    \end{split}
\end{equation*}
where \(\bQQ_{k-1}\trans \vX_{k+1}\) also remains independent of the quantities from Lines~\ref{line:bcggspipiro:Ydiag}--\ref{line:bcgspipiro:Rdiag} in the (\(k-1\))-st iteration.
Thus, the calculation of \(\vQ_k\trans \vX_{k+1}\) is essential to avoid the need for the second synchronization point.
Therefore, our aim is to employ an alternative method to compute \(\vQ_k\trans \vX_{k+1}\) such that the inner products related to large matrices can be moved to right after Line~\ref{line:bcggspipiro:YY-Omega} of the (\(k-1\))-st iteration and then the second synchronization point in Line~\ref{line:bcggspipiro:SS} can be removed.

Here we use a similar way to compute \(\vQ_k\trans \vX_{k+1}\) as shown in~\cite[Figure 3]{YTHBSE2020} and \(\BCGSIROAone\) introduced by~\cite{CLMO2024-ls, CLRT2022}.
Notice that \(\vQ_k\) satisfies
\[\vQ_k = (\vU_k - \bQQ_{k-1} \YY_{1:k-1,k}) Y_{kk}^{-1}.\]
Furthermore, we have
\begin{equation} \label{eq:QkXk+1}
    \begin{split}
        \vQ_k\trans \vX_{k+1} &= Y_{kk}\itrans (\underbrace{\vU_k\trans \vX_{k+1}}_{=: P_k} - \YY_{1:k-1,k}\trans \underbrace{\bQQ_{k-1}\trans \vX_{k+1}}_{=: \vZ_{k-1}}) \\
        &= Y_{kk}\itrans (P_k - \YY_{1:k-1,k}\trans \vZ_{k-1}),
    \end{split}
\end{equation}
where these intermediate quantities \(P_k\) and \(\vZ_{k-1}\) can be computed together with Line~\ref{line:bcggspipiro:YY-Omega} of the (\(k-1\))-st iteration in Algorithm~\ref{alg:BCGSPIPIRO}.
Then we derive \(\BCGSPIPIROone\), as outlined in Algorithm~\ref{alg:BCGSPIPIRO1S}, by using~\eqref{eq:QkXk+1} to compute \(\vQ_k\trans \vX_{k+1}\) instead of Line~\ref{line:bcggspipiro:SS} of Algorithm~\ref{alg:BCGSPIPIRO}.
Consequently, the synchronization point required by Line~\ref{line:bcggspipiro:SS} of Algorithm~\ref{alg:BCGSPIPIRO} can be removed, and this new variant includes only one synchronization point per iteration, located at Line~\ref{line:bcgspipiro1s:sync} in Algorithm~\ref{alg:BCGSPIPIRO1S}.

\begin{algorithm}[htbp!]
    \caption{$[\bQQ, \RR] = \BCGSPIPIRO(\bXX, \IOAnoarg)$ \label{alg:BCGSPIPIRO}}
    \begin{algorithmic}[1]
        \State{$[\vQ_1, R_{11}] = \IOA{\vX_1}$}
        \For{$k = 1, \ldots, p-1$}
            \State{$\bmat{\SS_{1:k,k+1} \\ T_{k+1}} = \bmat{\bQQ_{k} & \vX_{k+1}}\trans \vX_{k+1}$} \SComment{First synchronization} \label{line:bcggspipiro:SS}
            \State{$S_{k+1,k+1} = \chol\bigl( T_{k+1} - \SS_{1:k,k+1}\trans \SS_{1:k,k+1} \bigr)$} \label{line:bcgspipiro:Sdiag}
            \State{$\vU_{k+1} = (\vX_{k+1} - \bQQ_{k} \SS_{1:k,k+1}) S_{k+1,k+1}^\inv$} \label{line:bcggspipiro:Uk}
            \State{$\bmat{\YY_{1:k,k+1} \\ \Omega_{k+1}} = \bmat{\bQQ_{k} & \vU_{k+1}}\trans \vU_{k+1}$} \SComment{Second synchronization} \label{line:bcggspipiro:YY-Omega}
            \State{$Y_{k+1,k+1} = \chol\bigl( \Omega_{k+1} - \YY_{1:k,k+1}\trans \YY_{1:k,k+1} \bigr)$} \label{line:bcggspipiro:Ydiag}
            \State{$\vQ_{k+1} = (\vU_{k+1} - \bQQ_{k} \YY_{1:k,k+1}) Y_{k+1,k+1}^\inv$} \label{line:bcggspipiro:Qk}
            \State{$\RR_{1:k,k+1} = \SS_{1:k,k+1} + \YY_{1:k,k+1} S_{k+1,k+1}$}
            \State{$R_{k+1,k+1} = Y_{k+1,k+1} S_{k+1,k+1}$} \label{line:bcgspipiro:Rdiag}
        \EndFor
        \State \Return{$\bQQ = [\vQ_1, \ldots, \vQ_p]$, $\RR = (R_{ij})$}
    \end{algorithmic}
\end{algorithm}

\begin{algorithm}[!tb]
    \caption{$[\bQQ, \RR] = \BCGSPIPIROone(\bXX, \IOAnoarg)$ \label{alg:BCGSPIPIRO1S}}
    \begin{algorithmic}[1]
        \State{$[\vQ_1, R_{11}] = \IOA{\vX_1}$}
        \State $\bmat{S_{12} \\ T_2} = \bmat{\vQ_1& \vX_2}\trans \vX_2$
        \State \diff{$S_{22} = \chol (T_2 - S_{12}\trans S_{12})$} \label{line:bcgspipiro1s:S22}
        \State \diff{$\vU_2 = (\vX_2 - \vQ_1 S_{12}) S_{22}^{-1}$}
        \State {$\begin{bmatrix} \YY_{12} & \vZ_{1} \\ \Omega_2 & P_2 \end{bmatrix} = \begin{bmatrix}\vQ_{1} & \vU_2 \end{bmatrix}\trans \begin{bmatrix}\vU_2 & \vX_{3} \end{bmatrix}$} \SComment{Synchronization}
        \Statex \diff{$T_{3} = \vX_{3}\trans \vX_{3}$}
        \State $Y_{22} = \chol(\Omega_2 - \YY_{12}\trans \YY_{12})$
        \State $\vQ_2 = (\vU_2 - \vQ_{1} \YY_{12}) Y_{22}^{-1}$
        \State $\RR_{12} = \SS_{12} + \YY_{12} \diff{S_{22}}$
        \State $R_{22} = Y_{22} \diff{S_{22}}$
        \For{$k = 2, \ldots, p-1$}
            \State $\SS_{1:k,k+1} =
                \begin{bmatrix} \vZ_{k-1} \\
                Y_{kk}\itrans \left( P_k - \YY_{1:k-1,k}\trans \vZ_{k-1} \right) \end{bmatrix}$ \label{line:bcgspipiro1s:computeS}
            \State \diff{$S_{k+1, k+1} = \chol (T_{k+1} - \SS_{1:k,k+1}\trans \SS_{1:k,k+1})$} \label{line:bcgspipiro1s:chol-1}
            \State \diff{$\vU_{k+1} = (\vX_{k+1} - \bQQ_k \SS_{1:k,k+1}) S_{k+1, k+1}^{-1}$} \label{line:1s:Uk}
            \If{$k < p-1$}
                \State $\begin{bmatrix} \YY_{1:k,k+1} & \vZ_{k} \\ \Omega_{k+1} & P_{k+1} \end{bmatrix} =
                \begin{bmatrix}\bQQ_{k} & \vU_{k+1} \end{bmatrix}\trans \begin{bmatrix}\vU_{k+1} & \vX_{k+2} \end{bmatrix}$ \SComment{Synchronization} \label{line:bcgspipiro1s:sync}
                \Statex $\qquad\quad T_{k+2} = \vX_{k+2}\trans \vX_{k+2}$
            \Else
                \State $\begin{bmatrix} \YY_{1:k,k+1} \\ \Omega_{k+1}\end{bmatrix} =
                \begin{bmatrix}\bQQ_{k} & \vU_{k+1} \end{bmatrix}\trans \vU_{k+1}$ 
            \EndIf 
            \State $Y_{k+1,k+1} = \chol(\Omega_{k+1} - \YY_{1:k,k+1}\trans \YY_{1:k,k+1})$
            \State $\vQ_{k+1} = (\vU_{k+1} - \bQQ_{k} \YY_{1:k,k+1}) Y_{k+1,k+1}^{-1}$ \label{line:Qk}
            \State $\RR_{1:k,k+1} = \SS_{1:k,k+1} + \YY_{1:k,k+1} \diff{S_{k+1,k+1}}$
            \State $R_{k+1,k+1} = Y_{k+1,k+1} \diff{S_{k+1,k+1}}$
        \EndFor
        \State \Return{$\bQQ = [\vQ_1, \ldots, \vQ_p]$, $\RR = (R_{ij})$}
    \end{algorithmic}
\end{algorithm}

\subsection{Loss of orthogonality of \texttt{BCGSI+P-1s}}
\label{sec:stability}
For the following analysis, we assume that the quantities computed by Algorithm~\ref{alg:BCGSPIPIRO1S} in the \(k\)-th iteration satisfy
\begin{equation} \label{eq:1s:bQQk}
    \norm{I - \bQQhat_{k}\trans \bQQhat_k} \leq \omega_k, \quad
    \norm{\bQQhat_k} \leq \sqrt{1 + \omega_k} \leq \sqrt{2}
\end{equation}
with \(\omega_k \in (0, 1)\), and
\begin{align}
    \SShat_{1:k,k+1}			&= \bQQhat_{k}\trans \vX_{k+1} + \DeltavS_{k+1},
    \quad \norm{\DeltavS_{k+1}}	\leq \deltaQX^{(k)} \norm{\vX_{k+1}}; \label{eq:1s:vSk1} \\
    \Shat_{k+1,k+1}\trans\Shat_{k+1,k+1} 		&= \That_{k+1} - \SShat_{1:k,k+1}\trans \SShat_{1:k,k+1} + \DeltaFone_{k+1} + \DeltaCone_{k+1}, \label{eq:1s:Sk1k1_chol} \\
    \norm{\DeltaFone_{k+1}}			&\leq \deltaSS^{(k)} \norm{\vX_{k+1}}^2 \label{eq:1s:Sk1k1_chol_errbound}
    \mbox{ and } \norm{\DeltaCone_{k+1}} \leq \deltacholone^{(k)} \norm{\vX_{k+1}}^2; \\
    \vVhat_{k+1} 					&= \vX_{k+1} - \bQQhat_{k} \SShat_{1:k,k+1} + \DeltavV_{k+1}
    \quad \text{with}
    \label{eq:1s:Vk1} \\
     \norm{\DeltavV_{k+1}} &\leq \deltaQS^{(k)} \norm{\vX_{k+1}}; \\
    \vUhat_{k+1} \Shat_{k+1,k+1}        &= \vVhat_{k+1} + \DeltavGone_{k+1},
    \quad \norm{\DeltavGone_{k+1}} \leq \deltaU^{(k)} \norm{\vUhat_{k+1}} \norm{\Shat_{k+1,k+1}} \label{eq:1s:Uk1}; \\
    \That_{k+2}					&= \vX_{k+2}\trans \vX_{k+2} + \DeltaT_{k+2},
    \quad \norm{\DeltaT_{k+2}} \leq \deltaXX^{(k)} \norm{\vX_{k+2}}^2; \label{eq:1s:Tk1} \\
    \YYhat_{1:k,k+1}            &= \bQQhat_{k}\trans \vUhat_{k+1} + \DeltavY_{k+1},
    \quad \norm{\DeltavY_{k+1}} \leq \deltaQU^{(k)} \norm{\vUhat_{k+1}}; \label{eq:1s:vYk1} \\
    \Omhat_{k+1}                     &= \vUhat_{k+1}\trans \vUhat_{k+1} + \DeltaOmega_{k+1},
    \quad \norm{\DeltaOmega_{k+1}} \leq \deltaUU^{(k)} \norm{\vUhat_{k+1}}^2 \label{eq:1s:Omk1}; \\ 
    \Yhat_{k+1,k+1}\trans \Yhat_{k+1,k+1}     &= \Omhat_{k+1} - \YYhat_{1:k,k+1}\trans \YYhat_{1:k,k+1} + \DeltaFtwo_{k+1} + \DeltaCtwo_{k+1} \label{eq:1s:Yk1k1_chol}, \\
    \norm{\DeltaFtwo_{k+1}}         &\leq \deltaYY^{(k)} \norm{\vUhat_{k+1}}^2
    \mbox{ and }
    \norm{\DeltaCtwo_{k+1}} \leq \deltacholtwo^{(k)} \norm{\vUhat_{k+1}}^2; \\
    \vWhat_{k+1}                    &= \vUhat_{k+1} - \bQQhat_{k} \YYhat_{1:k,k+1} + \DeltavW_{k+1}
    \quad \text{with}
    \label{eq:1s:Wk1} \\
    \norm{\DeltavW_{k+1}} &\leq \deltaQY^{(k)} \norm{\vUhat_{k+1}}; \\
    \vQhat_{k+1} \Yhat_{k+1,k+1}         &= \vWhat_{k+1} + \DeltavGtwo_{k+1}, \quad \norm{\DeltavGtwo_{k+1}} \leq \deltaQ^{(k)} \norm{\vQhat_{k+1}} \norm{\Yhat_{k+1,k+1}}
    \label{eq:1s:Yk1k1}
\end{align}
with rounding-error bounds \(\delta_*^{(k)} \in (0, 1)\) except for \(\deltaQX^{(k)}\).
The superscript \(k\) of \(\delta_*^{(k)}\) denotes that these \(\delta_*^{(k)}\) are rounding-error bounds from the \(k\)-th iteration.

For the following proof, we first summarize the estimations of two intermediate variables \(\norm{\vUhat_{k+1}}\) and \(\sigmin\bigl((I - \bQQhat_{k} \bQQhat_{k}\trans)\vUhat_{k+1}\bigr)\) from~\cite[Lemma~2, Theorem~2, Equation~(70)]{CLMO2024-P} in the following lemma.

\begin{lemma} \label{lem:pipi+result}
    Fix $k \in \{1, \ldots, p-1\}$ and suppose that \eqref{eq:1s:bQQk}--\eqref{eq:1s:Tk1} are satisfied with \(\delta_*^{(k)} \in (0, 1)\).
    Furthermore suppose that $\bQQhat_{k}$ and $\RRhat_{k}$ satisfy the following:
    \begin{align}
         \bQQhat_{k} \RRhat_{k}
         &= \bXX_{k} + \DeltabDD_{k},\quad \norm{\DeltabDD_{k}} \leq \deltaX^{(k-1)} \norm{\bXX_{k}}, \label{eq:lem:PIPI+:res}
    \end{align}
    with $\deltaX^{(k-1)}$, $\omega_{k} \in (0,1)$.
    Assume that
    \begin{equation} \label{eq:lem:PIPI+:delUkappa}
        8 \, \deltaU^{(k)} \kappa^2(\bXX_{k+1}) \leq 1
    \end{equation}
    and
    \begin{equation}  \label{eq:lem:PIPI+:kappa}
        \Bigl(\deltaX^{(k-1)} + 2 \, \omega_{k} + 2 \, \deltaSS^{(k)} + 2 \, \deltaXX^{(k-1)} + 2 \, \deltacholone^{(k)} \\
         + 18 \, \deltaQX^{(k)} + 8 \, \deltaQS^{(k)}\Bigr) \kappa^2(\bXX_{k+1}) \leq \frac{1}{2}.
    \end{equation}
    Then it holds that
    \begin{align}
        \norm{\vUhat_{k+1}} & \leq 2. \label{eq:lem:PIPI+:Uk}
    \end{align}

    Furthermore, if~\eqref{eq:1s:vYk1}--\eqref{eq:1s:Yk1k1} are satisfied with \(\delta_*^{(k)} \in (0, 1)\) and
    \begin{equation} \label{eq:thm:PIPI+:LOO:assump}
        \begin{split}
            & 2\,\bigl(4 \, \omega_{k} + \deltaSS^{(k)} + \deltaXX^{(k-1)} + \deltacholone^{(k)} + 7 \, \deltaQS^{(k)}
            + 13 \, \deltaQX^{(k)} + 14 \, \deltaU^{(k)}  \bigr) \kappa^2(\bXX_{k+1}) \\
            & + 4\,\bigl(2 \, \omega_{k} + \deltaUU^{(k)} + \deltaYY^{(k)} + \deltacholtwo^{(k)} + 3 \, \deltaQU^{(k)}\bigr) \leq \frac{1}{2},
        \end{split}
    \end{equation}
    then it holds that
    \begin{equation} \label{eq:lem:PIPI+:sigminUk1}
        \begin{split}
            \sigmin^2((I - \bQQhat_{k} \bQQhat_{k}\trans)\vUhat_{k+1}) & \geq 1 - 2\,\bigl(4 \, \omega_{k} + \deltaSS^{(k)} + \deltaXX^{(k-1)} + \deltacholone^{(k)} + 13 \, \deltaQX^{(k)} \\
            & \quad + 7 \, \deltaQS^{(k)} + 14 \, \deltaU^{(k)} \bigr) \kappa^2(\bXX_{k+1}).
        \end{split}
    \end{equation}
\end{lemma}

\medskip
Using standard rounding-error analysis, deriving these bounds \(\delta_*\) is quite straightforward, except when it comes to \(\deltaQX\) in~\eqref{eq:1s:vSk1}.
The reason stems from \(\SShat_{1:k,k+1}\) being computed via Line~\ref{line:bcgspipiro1s:computeS} of Algorithm~\ref{alg:BCGSPIPIRO1S} rather than directly through \(\bQQhat_{k}\trans \vX_{k+1}\).
Thus, we also assume that
\begin{align}
    \SShat_{1:k-1, k+1} &= \vZhat_{k-1} = \bQQhat_{k-1}\trans \vX_{k+1} + \Delta \vZ_{k-1}, \qquad \norm{\Delta\vZ_{k-1}} \leq \deltaZ^{(k-1)} \norm{\vX_{k+1}}; \label{eq:1s:SSk-1k+1}\\
    \Phat_{k}                     &= \vUhat_{k}\trans \vX_{k+1} + \DeltaP_{k},
    \qquad \norm{\DeltaP_{k}} \leq \deltaP^{(k-1)} \norm{\vUhat_{k}} \norm{\vX_{k+1}} \label{eq:1s:Pk}; \\ 
    \vJhat_k &= \Phat_k - \YYhat_{1:k-1,k}\trans \vZhat_{k-1} + \DeltaJ_k, \qquad \norm{\DeltaJ_k} \leq \deltaJ^{(k)} \norm{\vUhat_{k}} \norm{\vX_{k+1}}; \label{eq:1s:Jk} \\
    \Yhat_{kk}\trans \SShat_{k, k+1} &= \vJhat_k + \Delta E_{k, k+1}, \qquad \norm{\Delta E_{k,k+1}} \leq \deltaS^{(k)} \norm{\Yhat_{kk}} \norm{\SShat_{k, k+1}}. \label{eq:1s:YkkSkk1} 
\end{align}
Lemma~\ref{lem:pipiro1s:SShat} that follows provides an estimation for \(\deltaQX\).

\begin{lemma} \label{lem:pipiro1s:SShat}
    Fix \(k \in \{2, \ldots, p-1\}\) and assume that~\eqref{eq:1s:bQQk}, \eqref{eq:1s:SSk-1k+1}--\eqref{eq:1s:YkkSkk1} hold with
    \begin{equation} \label{eq:lem-1s:assump}
        4\sqrt{2}\, \deltaS^{(k)} + 8 \,\deltaQ^{(k-1)} + 2\sqrt{2}\, \deltaQY^{(k-1)} + 2\sqrt{2}\, \deltaP^{(k-1)} + 2\sqrt{2}\, \deltaJ^{(k)} + 5\,\sqrt{2}\, \deltaZ^{(k-1)} \leq \frac{1}{2}.
    \end{equation}
    Furthermore, assume that~\eqref{eq:1s:vYk1}--\eqref{eq:1s:Yk1k1}, \eqref{eq:lem:PIPI+:Uk}, and~\eqref{eq:lem:PIPI+:sigminUk1} hold for the \((k-1)\)-st iteration with
    \begin{equation} \label{eq:lem-1s:assumpfork-1}
        \begin{split}
            & 2\,\bigl(4 \, \omega_{k-1} + \deltaSS^{(k-1)} + \deltaXX^{(k-2)} + \deltacholone^{(k-1)} + 13 \, \deltaQX^{(k-1)} + 7 \, \deltaQS^{(k-1)} + 14 \, \deltaU^{(k-1)} \bigr) \kappa^2(\bXX_{k}) \\
            & \quad
            + 4\, \bigl(\deltaUU^{(k-1)} + 2 \sqrt{2}\, \deltaQU^{(k-1)} + \deltaYY^{(k-1)} + \deltacholtwo^{(k-1)} + 2\, \omega_{k-1}\bigr) \leq \frac{1}{2}.
        \end{split}
    \end{equation}
    Then it holds that
    \begin{equation} \label{eq:lem:SShat}
        \begin{split}
            \SShat_{1:k, k+1} &= \bQQhat_{k}\trans \vX_{k+1} + \Delta \vS_{k+1}
        \end{split}
    \end{equation}
    with \(\Delta \vS_{k+1} := \bmat{\Delta \vZ_{k-1}& \Delta S_{k, k+1}}\) satisfying
    \[
        \norm{\Delta S_{k, k+1}} \leq \sqrt{2}\, \bigl(4\sqrt{2}\, \deltaQ^{(k-1)} + 2\, \deltaQY^{(k-1)} + 2\, \deltaP^{(k-1)} + 2\, \deltaJ^{(k)} + 5\, \deltaZ^{(k-1)} + 16\, \deltaS^{(k)}\bigr) \norm{\vX_{k+1}}.
    \]
\end{lemma}

\medskip
\begin{proof}
    From~\eqref{eq:1s:SSk-1k+1}, to prove~\eqref{eq:lem:SShat}, we only need to estimate \(\Delta S_{k, k+1}\) satisfying
    \begin{equation}
        \SShat_{k, k+1} = \vQhat_k\trans \vX_{k+1} + \Delta S_{k, k+1}.
    \end{equation}
    Applying \(\Yhat_{kk}\itrans\) to both sides of~\eqref{eq:1s:YkkSkk1}, we obtain
    \begin{equation} \label{eq:lemSShat:proof:SShatkk1-0}
        \SShat_{k, k+1} = \Yhat_{kk}\itrans \vJhat_k + \Yhat_{kk}\itrans \Delta E_{k, k+1}.
    \end{equation}
    Furthermore, substituting \(\Phat_k\) involved in~\eqref{eq:1s:Jk} by~\eqref{eq:1s:Pk} and then substituting \(\vJhat_k\) in~\eqref{eq:lemSShat:proof:SShatkk1-0} by~\eqref{eq:1s:Jk}, \(\Delta S_{k, k+1}\) satisfies
    \begin{equation*}
        \begin{split}
            \SShat_{k, k+1} &= \Yhat_{kk}\itrans \Phat_k - \Yhat_{kk}\itrans \YYhat_{1:k-1,k}\trans \vZhat_{k-1} + \Yhat_{kk}\itrans \DeltaJ_k + \Yhat_{kk}\itrans \Delta E_{k, k+1} \\
            &= \Yhat_{kk}\itrans (\vUhat_{k} - \bQQhat_{k-1} \YYhat_{1:k-1,k})\trans \vX_{k+1} + \Yhat_{kk}\itrans \DeltaP_{k} - \Yhat_{kk}\itrans \YYhat_{1:k-1,k}\trans \Delta \vZ_{k-1} \\
            &\quad + \Yhat_{kk}\itrans \DeltaJ_k + \Yhat_{kk}\itrans \Delta E_{k, k+1}.
        \end{split}
    \end{equation*}
    Using~\eqref{eq:1s:Wk1} to substitute \(\vWhat_k\) for \(\vUhat_{k} - \bQQhat_{k-1} \YYhat_{1:k-1,k}\), we obtain
    \begin{equation*}
        \begin{split}
            \SShat_{k, k+1} &= \Yhat_{kk}\itrans \vWhat_k \trans \vX_{k+1} - \Yhat_{kk}\itrans \DeltavW_k \trans \vX_{k+1} + \Yhat_{kk}\itrans \DeltaP_{k} - \Yhat_{kk}\itrans \YYhat_{1:k-1,k}\trans \Delta \vZ_{k-1} \\
            &\quad + \Yhat_{kk}\itrans \DeltaJ_k + \Yhat_{kk}\itrans \Delta E_{k, k+1}.
        \end{split}
    \end{equation*}
    Furthermore, using~\eqref{eq:1s:Yk1k1} to substitute \(\vQhat_k\trans\) for \(\Yhat_{kk}\itrans \vWhat_k \trans\), we have
    \begin{equation} \label{eq:lem-1s-proof:SShatkk1}
        \begin{split}
            \SShat_{k, k+1} &= \vQhat_k \trans \vX_{k+1} + \Delta S_{k, k+1},
        \end{split}
    \end{equation}
    where \(\Delta S_{k, k+1}\) is defined by
    \begin{equation}
        \begin{split}
            \Delta S_{k, k+1} &:= - \Yhat_{kk}\itrans (\DeltavGtwo_{k})\trans \vX_{k+1} - \Yhat_{kk}\itrans \DeltavW_k \trans \vX_{k+1} + \Yhat_{kk}\itrans \DeltaP_{k} \\
            &\quad - \Yhat_{kk}\itrans \YYhat_{1:k-1,k}\trans \Delta \vZ_{k-1} + \Yhat_{kk}\itrans \DeltaJ_k + \Yhat_{kk}\itrans \Delta E_{k, k+1}.
        \end{split}
    \end{equation}
    From~\eqref{eq:1s:Yk1k1}, \eqref{eq:1s:Wk1}, \eqref{eq:1s:Pk}--\eqref{eq:1s:YkkSkk1}, \eqref{eq:1s:SSk-1k+1}, \(\Delta S_{k, k+1}\) satisfies
    \begin{equation*}
        \begin{split}
            \norm{\Delta S_{k, k+1}} &\leq \deltaQ^{(k-1)} \norm{\vQhat_{k}} \norm{\Yhat_{kk}} \norm{\Yhat_{kk}^\inv} \norm{\vX_{k+1}} + \deltaZ^{(k-1)} \norm{\Yhat_{kk}^\inv} \norm{\YYhat_{1:k-1,k}} \norm{\vX_{k+1}} \\
            &\quad + \bigl(\deltaQY^{(k-1)} + \deltaP^{(k-1)} + \deltaJ^{(k)}\bigr) \norm{\Yhat_{kk}^\inv} \norm{\vUhat_{k}} \norm{\vX_{k+1}} \\
            &\quad + \deltaS^{(k)} \norm{\Yhat_{kk}^\inv} \norm{\Yhat_{kk}} \norm{\SShat_{k, k+1}}.
        \end{split}
    \end{equation*}
    Together with the bounds of \(\norm{\vQhat_{k}}\), \(\norm{\YYhat_{1:k-1,k}}\), and \(\norm{\vUhat_{k}}\) as shown in~\eqref{eq:1s:bQQk}, \eqref{eq:1s:vYk1}, and \eqref{eq:lem:PIPI+:Uk}, respectively, it follows that
    \begin{equation} \label{eq:lem-1s-proof-normDSkk1}
        \begin{split}
            \norm{\Delta S_{k, k+1}} &\leq \sqrt{2}\, \deltaQ^{(k-1)} \norm{\Yhat_{kk}} \norm{\Yhat_{kk}^\inv} \norm{\vX_{k+1}} + 5\, \deltaZ^{(k-1)} \norm{\Yhat_{kk}^\inv} \norm{\vX_{k+1}} \\
            &\quad + 2\, \bigl(\deltaQY^{(k-1)} + \deltaP^{(k-1)} + \deltaJ^{(k)}\bigr) \norm{\Yhat_{kk}^\inv} \norm{\vX_{k+1}} \\
            &\quad + \deltaS^{(k)} \norm{\Yhat_{kk}^\inv} \norm{\Yhat_{kk}} \norm{\SShat_{k, k+1}}.
        \end{split}
    \end{equation}
    Next, it remains to estimate \(\norm{\Yhat_{kk}}\), \(\norm{\Yhat_{kk}^\inv}\), and \(\norm{\SShat_{k, k+1}}\).

    First, we will treat \(\norm{\Yhat_{kk}}\) and \(\norm{\Yhat_{kk}^\inv}\).
    From~\eqref{eq:1s:vYk1}--\eqref{eq:1s:Yk1k1_chol}, it holds that
    \begin{equation} \label{eq:lem-1s-proof:YkkYkk}
        \begin{split}
            \Yhat_{kk}\trans \Yhat_{kk} &= \vUhat_k\trans (I - \bQQhat_{k-1} \bQQhat_{k-1}\trans) \vUhat_k + \Delta M_k \\
            &= \vUhat_k\trans (I - \bQQhat_{k-1} \bQQhat_{k-1}\trans) (I - \bQQhat_{k-1} \bQQhat_{k-1}\trans) \vUhat_k + \Delta \widetilde M_k,
        \end{split}
    \end{equation}
    where \(\Delta M_k := \Delta \Omega_k  + \vUhat_k\trans \bQQhat_{k-1} \DeltavY_k + \DeltavY_k\trans \bQQhat_{k-1}\trans \vUhat_k + \DeltaFtwo_k + \DeltaCtwo_k\) satisfies
    \begin{equation}
        \begin{split}
            \norm{\Delta M_k} &\leq \bigl(\deltaUU^{(k-1)} + 2 \sqrt{2}\, \deltaQU^{(k-1)}  + \deltaYY^{(k-1)} + \deltacholtwo^{(k-1)}\bigr) \norm{\vUhat_{k}}^2 \\
            &\leq 4\, \bigl(\deltaUU^{(k-1)} + 2 \sqrt{2}\, \deltaQU^{(k-1)}  + \deltaYY^{(k-1)} + \deltacholtwo^{(k-1)}\bigr)
        \end{split}
    \end{equation}
    and \(\Delta \widetilde M_k := \Delta M_k + \vUhat_k\trans \bQQhat_{k-1} \bQQhat_{k-1}\trans (I - \bQQhat_{k-1} \bQQhat_{k-1}\trans) \vUhat_k\) satisfies
    \begin{equation} \label{eq:lem-1s-proof:normtildeM}
        \norm{\Delta \widetilde M_k} \leq 4\, \bigl(\deltaUU^{(k-1)} + 2 \sqrt{2}\, \deltaQU^{(k-1)}  + \deltaYY^{(k-1)} + \deltacholtwo^{(k-1)} + 2\, \omega_{k-1}\bigr).
    \end{equation}
    This implies that
    \begin{equation} \label{eq:lem-1s-proof:normYkk}
        \norm{\Yhat_{kk}} \leq 2 \sqrt{3 + \deltaUU^{(k-1)} + 2 \sqrt{2}\, \deltaQU^{(k-1)}  + \deltaYY^{(k-1)} + \deltacholtwo^{(k-1)}} \leq 4.
    \end{equation}
    Combined~\eqref{eq:lem-1s-proof:YkkYkk}, \eqref{eq:lem-1s-proof:normtildeM}, and~\eqref{eq:lem:PIPI+:sigminUk1} with~\eqref{eq:lem:PIPI+:Uk} and the assumption~\eqref{eq:lem-1s:assumpfork-1}, it follows that
    \begin{equation} \label{eq:lem-1s-proof:sigminYkk}
        \begin{split}
            \sigmin^2(\Yhat_{kk}) &\geq \sigmin^2((I - \bQQhat_{k-1} \bQQhat_{k-1}\trans)\vUhat_k) - 4\, \bigl(\deltaUU^{(k-1)} + 2 \sqrt{2}\, \deltaQU^{(k-1)} \\
            & \quad + \deltaYY^{(k-1)} + \deltacholtwo^{(k-1)} + 2\, \omega_{k-1}\bigr) \\
            &\geq 1 - 2\,\bigl(4 \, \omega_{k-1} + \deltaSS^{(k-1)} + \deltaXX^{(k-2)} + \deltacholone^{(k-1)} + 13 \, \deltaQX^{(k-1)} \\
            & \quad + 7 \, \deltaQS^{(k-1)} + 14 \, \deltaU^{(k-1)} \bigr) \kappa^2(\bXX_{k})
            - 4\, \bigl(\deltaUU^{(k-1)} + 2 \sqrt{2}\, \deltaQU^{(k-1)} \\
            & \quad + \deltaYY^{(k-1)} + \deltacholtwo^{(k-1)} + 2\, \omega_{k-1}\bigr) \\
            & \geq \frac{1}{2},
        \end{split}
    \end{equation}
    where the first inequality is derived from~\eqref{eq:lem-1s-proof:YkkYkk} and~\eqref{eq:lem-1s-proof:normtildeM}, the second inequality is derived from~\eqref{eq:lem:PIPI+:sigminUk1}, and the last one is from the assumption~\eqref{eq:lem-1s:assumpfork-1}.
    Then \eqref{eq:lem-1s-proof:sigminYkk} gives that
    \begin{equation} \label{eq:lem-1s-proof:normYkkinv}
        \norm{\Yhat_{kk}^\inv} = \frac{1}{\sigmin(\Yhat_{kk})} \leq \sqrt{2}.
    \end{equation}

    Then we will bound \(\norm{\SShat_{k, k+1}}\).
    From~\eqref{eq:lem-1s-proof:SShatkk1}, \eqref{eq:lem-1s-proof-normDSkk1}, \eqref{eq:lem-1s-proof:normYkk}, and~\eqref{eq:lem-1s-proof:normYkkinv}, we obtain
    \begin{equation} \label{eq:lem-1s-proof:normDSkk1-final}
        \begin{split}
            \norm{\Delta S_{k, k+1}} &\leq 8\, \deltaQ^{(k-1)} \norm{\vX_{k+1}} + 2\sqrt{2}\, \bigl(\deltaQY^{(k-1)} + \deltaP^{(k-1)} + \deltaJ^{(k)}\bigr) \norm{\vX_{k+1}} \\
            &\quad + 5\sqrt{2}\, \deltaZ^{(k-1)} \norm{\vX_{k+1}} + 4\sqrt{2}\, \deltaS^{(k)} \norm{\SShat_{k, k+1}},
        \end{split}
    \end{equation}
    and then by the assumption~\eqref{eq:lem-1s:assump},
    \begin{equation}
        \begin{split}
            \norm{\SShat_{k, k+1}} &\leq \frac{\sqrt{2}\, \bigl(1 + 4 \sqrt{2}\, \deltaQ^{(k-1)} + 2\, \deltaQY^{(k-1)} + 2\, \deltaP^{(k-1)} + 2\, \deltaJ^{(k)} + 5\, \deltaZ^{(k-1)}\bigr)}{1 - 4\sqrt{2}\, \deltaS^{(k)}} \norm{\vX_{k+1}} \\
            &\leq 4 \norm{\vX_{k+1}},
        \end{split}
    \end{equation}
    which concludes the proof combined with~\eqref{eq:lem-1s-proof:normDSkk1-final}.
\end{proof}

As shown in Section~\ref{sec:algo-onesync}, the sole distinction between Algorithm~\ref{alg:BCGSPIPIRO1S} and Algorithm~\ref{alg:BCGSPIPIRO}, i.e., \(\BCGSPIPIRO\), is the manner in which the former computes \(\SS_{1:k, k+1}\) by Line~\ref{line:bcgspipiro1s:computeS} of Algorithm~\ref{alg:BCGSPIPIRO1S}.
Combining Lemma~\ref{lem:pipiro1s:SShat} with Lemma~\ref{lem:pipi+result}, we have proven that~\eqref{eq:1s:vSk1} holds with \(\bigO(\macheps)\) in the \(k\)-th iteration, if \(\bigO(\macheps) \kappa^2(\bXX_{k+1}) < 1\) and~\eqref{eq:1s:vSk1}--\eqref{eq:1s:Yk1k1} hold with \(\delta_*^{(k-1)} = \bigO(\macheps)\) in the \((k-1)\)-st iteration.
Notice that the first iteration of Algorithm~\ref{alg:BCGSPIPIRO1S} is as same as the first iteration of \BCGSPIPIRO.
This means that~\eqref{eq:1s:vSk1}--\eqref{eq:1s:Yk1k1} hold with \(\delta_*^{(1)} = \bigO(\macheps)\) in the first iteration of Algorithm~\ref{alg:BCGSPIPIRO1S} from standard rounding-error analysis.
Then together with Lemma~\ref{lem:pipiro1s:SShat}, \eqref{eq:1s:vSk1}--\eqref{eq:1s:Yk1k1} also hold with \(\delta_*^{(k)} = \bigO(\macheps)\) in the \(k\)-th iteration of Algorithm~\ref{alg:BCGSPIPIRO1S} by induction.
Therefore, the LOO analysis presented in~\cite[Corollaries 3 and 4]{CLMO2024-P}, whose proof is also based on~\eqref{eq:1s:vSk1}--\eqref{eq:1s:Yk1k1}, gives the LOO result for Algorithm~\ref{alg:BCGSPIPIRO1S} as follows.

\begin{theorem} \label{thm:pipiro1s}
    Assume that \(\bQQhat_k\) and \(\RRhat_k\) are computed by Algorithm~\ref{alg:BCGSPIPIRO1S} with \(\bigO(\macheps) \kappa^2(\bXX_k) \leq 1/2\).
    Assume also that for all $\vX \in \mathcal{R}^{m \times s}$ with $\kappa(\vX) \leq \kappa(\bXX)$, $[\vQhat, \Rhat] = \IOA{\vX}$ satisfy
    \begin{align*}
        & \Rhat^T \Rhat = \vX^T \vX + \DeltavD, \quad \norm{\DeltavD} \leq \bigO(\macheps) \norm{\vX}^2 \\ 
        & \vX + \Delta \vX = \vQhat \Rhat , \quad \norm{\Delta \vX} \leq \bigO(\macheps) \norm{\vX} \mbox{ and} \\
        & \norm{I - \vQhat^T \vQhat} \leq \frac{\bigO(\macheps)}{1 - \bigO(\macheps) \kappa^2(\vX)}. 
    \end{align*} 
    Then for all $k \in \{1, \ldots, p\}$, \(\bQQhat_k\) and \(\RRhat_k\) satisfy
    \begin{align}
        &\bXX_{k} + \Delta\bXX_{k} = \bQQhat_{k} \RRhat_{k}, \quad \norm{\Delta\bXX_{k}} \leq \bigO(\macheps) \norm{\bXX_{k}}, \\
        &\norm{I - \bQQhat_k\trans \bQQhat_k} \leq \bigO(\macheps).  \label{eq:thm-pipiro1s:loo}
    \end{align}
\end{theorem}

\section{Two-sync Pythagorean-based BCGS}
\label{sec:algo-twosync}
In this section, we propose a two-sync reorthogonalized BCGS, as summarized in Algorithm~\ref{alg:BCGSPIPIRO2S}, that shares the same LOO property as \BCGSIROA described in~\cite{CLMO2024-ls}, which is equivalent to \BCGSIRO when using \HouseQR as \(\IOAnoarg\).
By adding one additional synchronization, we will effectively relax the constraint on $\kappa(\bXX_{k+1})$. 

\subsection{Algorithm}
Revisiting the analysis of Algorithm~\ref{alg:BCGSPIPIRO1S}, we find that the computation of \(S_{k+1,k+1}\) and \(\vU_{k+1}\), as described in Line~\ref{line:bcgspipiro1s:chol-1}, involves \(\vX_{k+1}\trans \vX_{k+1}\), which introduces the constraint \(\bigO(\macheps) \kappa^2(\bXX_{k+1}) \leq 1/2\) on the LOO bound for Algorithm~\ref{alg:BCGSPIPIRO1S}.
To address this, we adopt a more stable \(\IOonenoarg\), such as \(\HouseQR\) and \(\TSQR\), in Line~\ref{line:2s:io1} of Algorithm~\ref{alg:BCGSPIPIRO2S}.
This replaces the Pythagorean-based Cholesky QR method used in Lines~\ref{line:bcgspipiro1s:chol-1} and~\ref{line:1s:Uk} of Algorithm~\ref{alg:BCGSPIPIRO1S}.
This adjustment aims to circumvent operations involving \(\vX_{k+1}\trans \vX_{k+1}\) during the \(k\)-th iteration and achieve \(\bigO(\macheps)\) LOO under the condition \(\bigO(\macheps) \kappa(\bXX_{k+1}) \leq 1/2\) instead of \(\bigO(\macheps) \kappa^2(\bXX_{k+1}) \leq 1/2\), albeit at the cost of adding an additional synchronization point in each iteration.
We call the resulting algorithm \BCGSPIPIROtwo. 

\begin{algorithm}[htb!]
    \caption{$[\bQQ, \RR] = \BCGSPIPIROtwo(\bXX, \IOAnoarg, \IOonenoarg)$ \label{alg:BCGSPIPIRO2S}}
    \begin{algorithmic}[1]
        \State{$[\vQ_1, R_{11}] = \IOA{\vX_1}$}
        \State $S_{12} = \vQ_1\trans \vX_2$
        \State \diff{$[\vU_2, S_{22}] = \IOone{\vX_2 - \vQ_1 S_{12}}$} \SComment{First synchronization}
        \State {$\begin{bmatrix} \YY_{12} & \vZ_{1} \\ \Omega_2 & P_2 \end{bmatrix} = \begin{bmatrix}\vQ_{1} & \vU_2 \end{bmatrix}\trans \begin{bmatrix}\vU_2 & \vX_{3} \end{bmatrix}$} \SComment{Second synchronization}
        \State $Y_{22} = \chol(\Omega_2 - \YY_{12}\trans \YY_{12})$
        \State $\vQ_2 = (\vU_2 - \vQ_{1} \YY_{12}) Y_{22}^{-1}$
        \State $\RR_{12} = \SS_{12} + \YY_{12} \diff{S_{22}}$
        \State $R_{22} = Y_{22} \diff{S_{22}}$
        \For{$k = 2, \ldots, p-1$}
            \State $\SS_{1:k,k+1} =
                \begin{bmatrix} \vZ_{k-1} \\
                Y_{kk}\itrans \left( P_k - \YY_{1:k-1,k}\trans \vZ_{k-1} \right) \end{bmatrix}$ \label{line:bcgspipiro2s:computeS}
            \State $[\vU_{k+1}, S_{k+1, k+1}] = \IOone{\vX_{k+1} - \bQQ_k \SS_{1:k,k+1}}$ \SComment{First synchronization} \label{line:2s:io1}
            \If{$k < p-1$}
                \State $\begin{bmatrix} \YY_{1:k,k+1} & \vZ_{k} \\ \Omega_{k+1} & P_{k+1} \end{bmatrix} =
                \begin{bmatrix}\bQQ_{k} & \vU_{k+1} \end{bmatrix}\trans \begin{bmatrix}\vU_{k+1} & \vX_{k+2} \end{bmatrix}$ \SComment{Second synchronization} \label{line:bcgspipiro2s:sync}
            \Else
                \State $\begin{bmatrix} \YY_{1:k,k+1} \\ \Omega_{k+1} \end{bmatrix} =
                \begin{bmatrix}\bQQ_{k} & \vU_{k+1} \end{bmatrix}\trans \vU_{k+1}$
            \EndIf
            \State $Y_{k+1,k+1} = \chol(\Omega_{k+1} - \YY_{1:k,k+1}\trans \YY_{1:k,k+1})$
            \State $\vQ_{k+1} = (\vU_{k+1} - \bQQ_{k} \YY_{1:k,k+1}) Y_{k+1,k+1}^{-1}$ \label{line:bcgspipiro2s:Qk}
            \State $\RR_{1:k,k+1} = \SS_{1:k,k+1} + \YY_{1:k,k+1} \diff{S_{k+1,k+1}}$
            \State $R_{k+1,k+1} = Y_{k+1,k+1} \diff{S_{k+1,k+1}}$
        \EndFor
        \State \Return{$\bQQ = [\vQ_1, \ldots, \vQ_p]$, $\RR = (R_{ij})$}
    \end{algorithmic}
\end{algorithm}

\subsection{Loss of orthogonality of \texttt{BCGSI+P-2S}}
Comparing Algorithm~\ref{alg:BCGSPIPIRO2S}, i.e., \BCGSPIPIROtwo, with Algorithm~\ref{alg:BCGSPIPIRO1S}, the only difference is that we can employ a more stable \(\IOonenoarg\) in Line~\ref{line:2s:io1} of Algorithm~\ref{alg:BCGSPIPIRO2S} instead of using Pythagorean-based Cholesky QR in Lines~\ref{line:bcgspipiro1s:chol-1} and~\ref{line:1s:Uk} of Algorithm~\ref{alg:BCGSPIPIRO1S}.
Thus, in the \(k\)-th iteration, we can also assume that the quantities computed by Algorithm~\ref{alg:BCGSPIPIRO2S} satisfy~\eqref{eq:1s:bQQk}, \eqref{eq:1s:vSk1}, \eqref{eq:1s:Vk1}, \eqref{eq:1s:vYk1}--\eqref{eq:1s:Yk1k1}, and~\eqref{eq:1s:SSk-1k+1}--\eqref{eq:1s:YkkSkk1}. 
For \(\IOonenoarg\), we assume that
\begin{equation} \label{eq:2s:IO2-Vk1}
    \begin{split}
        & \vVhat_{k+1} + \DeltavVhat_{k+1} = \vUhat_{k+1} \SShat_{k+1,k+1}, \quad \norm{\DeltavVhat_{k+1}} \leq \deltaU^{(k)} \norm{\vVhat_{k+1}}, \\
        & \norm{I - \vUhat_{k+1}\trans \vUhat_{k+1}} \leq \omega_{qr}, \quad
        \norm{\vUhat_{k+1}} \leq \sqrt{1 + \omega_{qr}}, \quad
        \sigmin(\vUhat_{k+1}) \geq \sqrt{1 - \omega_{qr}}.
    \end{split}
\end{equation}
Similarly as in the proof of \BCGSPIPIROone, i.e., Lemma~\ref{lem:pipi+result}, we first need to bound \(\vUhat_{k+1}\) and \(\sigmin\bigl((I - \bQQhat_{k} \bQQhat_{k}\trans)\vUhat_{k+1}\bigr)\).

\begin{lemma} \label{lem:2s:Uk1}
    Fix $k \in \{1, \ldots, p-1\}$ and assume that for all $\vX \in \mathbb{R}^{m \times s}$ with $\bigO(\macheps) \kappa(\vX) \leq 1/2$, the following hold for $[\vQhat, \Rhat] = \IOone{\vX}$:
    \begin{align}
        & \vX + \Delta \vX = \vQhat \Rhat,
        \quad \norm{\Delta \vX} \leq \bigO(\macheps) \norm{\vX}, \notag \\
        & \norm{I - \vQhat^T \vQhat} \leq \bigO(\macheps) \kappa^{\alpha_1}(\vX) \label{eq:lem:assump-IOone}
    \end{align}
    with \(\alpha_1 \leq 1\).
    Then it holds that
    \begin{align}
        \norm{\vUhat_{k+1}} & \leq \sqrt{2} \quad\text{and}\quad
        \sigmin(\vUhat_{k+1}) \geq \frac{\sqrt{2}}{2}. \label{eq:lem:2s:Uk1:Uk}
    \end{align}

    Furthermore, assume that~\eqref{eq:1s:bQQk} and~\eqref{eq:1s:vSk1} hold with \(\omega_k = \bigO(\macheps)\) and \(\deltaQX^{(k)} = \bigO(\macheps)\) in the \(k\)-th iteration.
    If
    \begin{equation}
        \begin{split}
            & \bXX_{k} + \Delta \bXX_{k} = \bQQhat_{k} \RRhat_{k},
            \quad \norm{\Delta \bXX_{k}}
            \leq \bigO(\macheps) \norm{\bXX_{k}}
        \end{split}
    \end{equation}
    as well as
    \begin{equation}\label{eq:lem-bcgsiro-property:assump}
        \bigO(\macheps) \kappa(\bXX_{k+1}) + \bigO(\macheps) \leq \frac{1}{2}.
    \end{equation}
    Then it holds that
    \begin{equation} \label{eq:2s:Uk1:sigminUk1}
        \begin{split}
            \sigmin\bigl((I - \bQQhat_{k} \bQQhat_{k}\trans)\vUhat_{k+1}\bigr) & \geq \frac{\sqrt{2}}{2}
            - \bigO(\macheps) \kappa(\bXX_{k+1})
        \end{split}
    \end{equation}
    and
    \begin{equation} \label{eq:2s:loo}
        \norm{I - \bQQhat_{k+1}\trans \bQQhat_{k+1}} \leq \bigO(\macheps).
    \end{equation}
\end{lemma}

\medskip
\begin{proof}
     The assumption on \(\IOonenoarg\) means that \(\omega_{qr} \leq 1/2\) in~\eqref{eq:2s:IO2-Vk1}, which directly gives~\eqref{eq:lem:2s:Uk1:Uk}.
     From~\cite[Lemma~6 and Equation (58)]{CLMO2024-ls}, we have
     \begin{equation}
         \sigmin\bigl((I - \bQQhat_{k} \bQQhat_{k}\trans)\vUhat_{k+1}\bigr) \geq \sigmin (\vUhat_{k+1})
         - 2 \left(\omega_k + \bigO(\macheps) \right)(1 + \omega_k)^2 \kappa(\bXX_{k+1}),
     \end{equation}
     which proves~\eqref{eq:2s:Uk1:sigminUk1} combined with~\eqref{eq:lem:2s:Uk1:Uk}.
     
     It remains to prove~\eqref{eq:2s:loo}, which satisfies
     \begin{equation}
         \norm{I - \bQQhat_{k+1}\trans \bQQhat_{k+1}} \leq \omega_k + 2 \norm{\bQQhat_{k}\trans \vQhat_{k+1}} + \norm{I - \vQhat_{k+1}\trans \vQhat_{k+1}}.
     \end{equation}
     Thus, we need to treat \(\norm{\bQQhat_{k}\trans \vQhat_{k+1}}\) and \(\norm{I - \vQhat_{k+1}\trans \vQhat_{k+1}}\), respectively.
     Notice that \(\delta_*^{(k)} = \bigO(\macheps)\) in~\eqref{eq:1s:vSk1}, \eqref{eq:1s:Vk1}, and \eqref{eq:1s:vYk1}--\eqref{eq:1s:Yk1k1}, according to standard rounding-error analysis.
     The proof of these two norms is similar to the proof of~\cite[Theorem 2]{CLMO2024-P}  (see~\cite[Equation (74) and (78)]{CLMO2024-P} for details), i.e.,
     \begin{align*}
         & \norm{\bQQhat_{k}\trans \vQhat_{k+1}} \leq \bigO(\macheps), \quad\text{and}\quad
         \norm{I - \vQhat_{k+1}\trans \vQhat_{k+1}} \leq \bigO(\macheps),
     \end{align*}
     which concludes the proof of~\eqref{eq:2s:loo}.
\end{proof}

Lemma~\ref{lem:2s:Uk1} outlines the condition for \(\IOonenoarg\), requiring that the computed \(Q\)-factor of \(\IOonenoarg\) satisfies~\eqref{eq:lem:assump-IOone} with \(\alpha_1 \leq 1\).
Table~\ref{table:IOone-twosync} presents the provable LOO bounds, represented as \(\norm{I - \vQhat^T \vQhat} \leq \bigO(\macheps) \kappa^{\alpha_1}(\vX)\), which depends on \(\alpha_1\) corresponding to different choices of \(\IOonenoarg\).
Using \HouseQR, \TSQR, or \MGS as \(\IOonenoarg\) meets the condition set in Lemma~\ref{lem:2s:Uk1}, namely, \(\alpha_1 \leq 1\).

\begin{table}[!htbp]
    \centering
    \begin{tabular}{cc} \hline
        $\IOonenoarg$     & $\alpha_1$  \\ \hline
        \HouseQR          & 0           \\
        \TSQR             & 0           \\
        \MGS              & 1           \\
        \CholQR           & 2           \\ \hline
    \end{tabular}
    \caption{Values of \(\alpha_1\) for common \(\IOonenoarg\) choices. \label{table:IOone-twosync}}
\end{table}

Comparing Lemma~\ref{lem:2s:Uk1} with Lemma~\ref{lem:pipi+result}, these two lemmas both bound \(\vUhat_{k+1}\) by a small constant, but Lemma~\ref{lem:2s:Uk1} only requires \(\bigO(\macheps) \kappa(\bXX_{k+1}) + \bigO(\macheps) \leq 1/2\) instead of \(\bigO(\macheps) \kappa^2(\bXX_{k+1}) \leq 1/2\).
We are now prepared to give the following theorem to describe the LOO of \BCGSPIPIROtwo.

\begin{theorem} \label{thm:pipiro2s}
    Assume that \(\bQQhat_k\) and \(\RRhat_k\) are computed by Algorithm~\ref{alg:BCGSPIPIRO2S} with \(\bigO(\macheps) \kappa(\bXX_k) \leq 1/2\).
    Assume also that for $\bXX_1$ with $\kappa(\bXX_1) \leq \kappa(\bXX)$, $[\vQhat, \Rhat] = \IOA{\bXX_1}$ satisfies
    \begin{align*}
        & \bXX_1 + \Delta \bXX_1 = \vQhat \Rhat,
            \quad \norm{\Delta \bXX_1} \leq \bigO(\macheps) \norm{\bXX_1}, \\
            & \norm{I - \vQhat^T \vQhat} \leq \bigO(\macheps), 
    \end{align*} 
    and for all $\vX \in \mathcal{R}^{m \times s}$ with $\kappa(\vX) \leq \kappa(\bXX)$, $[\vQhat, \Rhat] = \IOone{\vX}$ satisfies
    \begin{equation*}
        \begin{split}
            & \vX + \Delta \vX = \vQhat \Rhat,
            \quad \norm{\Delta \vX} \leq \bigO(\macheps) \norm{\vX}, \\
            & \norm{I - \vQhat^T \vQhat} \leq \bigO(\macheps) \kappa^{\alpha_1}(\vX)
        \end{split}
    \end{equation*}
    with \(\alpha_1 \leq 1\).
    Then for all $k \in \{1, \ldots, p\}$, \(\bQQhat_k\) and \(\RRhat_k\) satisfy
    \begin{align}
        &\bXX_{k} + \Delta\bXX_{k} = \bQQhat_{k} \RRhat_{k}, \quad \norm{\Delta\bXX_{k}} \leq \bigO(\macheps) \norm{\bXX_{k}}, \label{eq:thm-pipiro2s:res} \\
        &\norm{I - \bQQhat_k\trans \bQQhat_k} \leq \bigO(\macheps). \label{eq:thm-pipiro2s:loo}
    \end{align}
\end{theorem}

\medskip
\begin{proof}
    We will prove~\eqref{eq:thm-pipiro2s:res} and~\eqref{eq:thm-pipiro2s:loo} by induction.
    For the base case \(k = 1\), it is clear that \([\vQhat_{1}, \Rhat_{11}] = \IOA{\bXX_1}\) satisfies~\eqref{eq:thm-pipiro2s:res} and~\eqref{eq:thm-pipiro2s:loo} due to the assumption on \(\IOAnoarg\).
    In addition, \eqref{eq:1s:vSk1} holds with \(\deltaQX^{(1)} = \bigO(\macheps)\) from standard rounding-error analysis.
    Also, note that \(\delta_*^{(k)} = \bigO(\macheps)\) in~\eqref{eq:1s:vSk1}, \eqref{eq:1s:Vk1}, and \eqref{eq:1s:vYk1}--\eqref{eq:1s:Yk1k1}, according to standard rounding-error analysis.
    Then we assume that~\eqref{eq:thm-pipiro2s:res} and~\eqref{eq:thm-pipiro2s:loo} hold, and also \eqref{eq:1s:vSk1} holds with \(\deltaQX^{(k)} = \bigO(\macheps)\), for the case \(k = i-1\).
    Our aim is to prove that these also hold for \(k = i\). 

    From Lemma~\ref{lem:2s:Uk1} and the induction hypothesis, we have
    \begin{align}
        \norm{\vUhat_i} &\leq \sqrt{2}, \label{eq:thm-2s-proof:normUi} \\
        \sigmin\bigl((I - \bQQhat_{i-1} \bQQhat_{i-1}\trans)\vUhat_{i}\bigr) & \geq \frac{\sqrt{2}}{2}
        - \bigO(\macheps) \kappa(\bXX_{i}), \\
        \norm{I - \bQQhat_{i}\trans \bQQhat_{i}} &\leq \bigO(\macheps). \label{eq:thm-2s-proof:looi}
    \end{align}
    Together with Lemma~\ref{lem:pipiro1s:SShat}, it follows that~\eqref{eq:1s:vSk1} holds with \(\deltaQX^{(k)} = \bigO(\macheps)\) for \(k = i\).
    Based on~\eqref{eq:thm-2s-proof:normUi}, \eqref{eq:thm-2s-proof:looi}, and~\eqref{eq:1s:vSk1} holding with \(\deltaQX^{(i)} = \bigO(\macheps)\), using~\cite[Corollary 4]{CLMO2024-P} we can obtain that~\eqref{eq:thm-pipiro2s:res} also holds for \(k = i\). 
\end{proof}
\section{Application of low-sync BCGS in \(s\)-step GMRES}
\label{sec:orthoGMRES}
An essential use of low-synchronization BCGS variants is in communication-avoiding Krylov subspaces methods.
We take communication-avoiding GMRES, also known as \(s\)-step GMRES, as an example.
GMRES is a popular Krylov subspace algorithm for solving nonsymmetric, nonsingular linear systems \(Ax = b\), which chooses \(x^{(k)} \in x^{(0)} + \krylov_k(A, r)\) to minimize \(\norm{A x^{(k)} - b}\) in the \(k\)-th iteration with a Krylov subspace \(\krylov_k(A, r) = \Span\{r, A r, \dotsc, A^{k-1} r\}\) and the initial residual \(r = b - A x^{(0)}\).

In each iteration, the \(s\)-step GMRES algorithm constructs a block Krylov basis
\begin{equation}
    \bmat{p_0(A) v& p_1(A) v& \dotsc & p_{s-1}(A) v}
\end{equation}
by given polynomials \(p_0\), \(p_1\), \(\dotsc\), \(p_{s-1}\), and employs a block orthogonalization method to orthogonalize them all at once.
We use \(B_{1:ks}\) to denote the Krylov basis after \(k\) iterations, and apply the block orthogonalization method to \(\bXX_k = \bmat{r & M_L^{-1}AM_R^{-1}B_{1:ks}}\), where \(M_L\) and \(M_R\) are the left and right preconditioners, respectively.
As demonstrated in~\cite{CM2024}, the orthogonalization method employed in \(s\)-step GMRES necessitates the following properties:
\begin{enumerate}
    \item Column-wise backward stability: if \(\bigO(\macheps) \kappa(\bXX) \leq 1\), the computed results satisfy
    \begin{equation*}
        \bXX + \Delta\bXX = \bQQhat \RRhat, \qquad \norm{\Delta X_i} \leq \bigO(\macheps) \norm{X_i} \quad \forall\, i = 1, 2, \dotsc, p
    \end{equation*}
    or
    \begin{equation*}
        \bXX + \Delta\bXX = \bQQtil \RRhat, \qquad \norm{\Delta X_i} \leq \bigO(\macheps) \norm{X_i} \quad \forall\, i = 1, 2, \dotsc, p
    \end{equation*}
    with an exactly orthonormal matrix \(\bQQtil\).
    \item Loss of orthogonality: if \(\bigO(\macheps) \kappa^\alpha(\bXX) \leq 1\), the computed \(Q\)-factor satisfies
    \begin{equation} \label{eq:loo}
        \norm{I - \bQQhat\trans \bQQhat} \leq \bigO(\macheps) \kappa^\beta(\bXX).
    \end{equation}
\end{enumerate}
Note that the first property is relatively straightforward to verify and is satisfied by all popular orthogonalization methods, including \(\BCGSIROA\) shown in~\cite[Algorithm~2]{CLMO2024-ls}, \(\BCGSPIPIROtwo\), and also \(\BCGSPIPIROone\).
In the case of the second property, it has been demonstrated in~\cite[Section 5.2]{CM2024} that the relative backward error can only be expected to be bounded by \(\bigO(\sqrt{\macheps})\) if \(\alpha = 2\) regardless of \(\beta\).
This makes orthogonalization methods with \(\alpha = 2\), such as \(\BCGSPIPIRO\) and \(\BCGSPIPIROone\), useless if one requires a relative backward error at the level \(\bigO(\macheps)\).

\subsection{An adaptive approach}
Here we propose an adaptive technique that combines both \(\BCGSPIPIROone\) and \(\BCGSPIPIROtwo\) to achieve \(\alpha = 1\) and \(\beta = 0\) in~\eqref{eq:loo} while requiring as few synchronization points as possible.
From Theorem~\ref{thm:pipiro2s}, \(\BCGSPIPIROtwo\) is adaptable to any choice of the first intraorthogonalization \(\IOAnoarg\) with \(\bigO(\macheps)\) LOO, while maintaining the same order of LOO.
Note that \BCGSPIPIROone has been proven to have \(\bigO(\macheps)\) LOO when the input matrix \(\bXX\) satisfies \(\bigO(\macheps) \kappa^2(\bXX) \leq 1/2\) according to Theorem~\ref{thm:pipiro1s}.
Consequently, it is natural to consider employing \BCGSPIPIROone to compute the orthogonalization of the first few well-conditioned block vectors.
Our approach initially utilizes \(\BCGSPIPIROone\), and subsequently switches to \(\BCGSPIPIROtwo\) for the remainder of the GMRES iterations once the condition number of the input matrix grows too large, ensuring \(\bigO(\macheps)\) LOO is maintained.

For determining the switching condition, it is natural to check if \(\bigO(\macheps) \kappa^2(\bXX_k) \leq 1/2\).
However, empirical observations reveal that LOO can still reach the level \(\bigO(\macheps)\) even when \(\bigO(\macheps) \kappa^2(\bXX_k) \leq 1/2\) is not satisfied.
Note that \(\bigO(\macheps) \kappa^2(\bXX_k) \leq 1/2\) is employed to ensure that \(\vU_k\) is well-conditioned in the LOO analysis of \BCGSPIPIROone.
Furthermore, the LOO can be guaranteed to be \(\bigO(\macheps)\) as long as \(\vU_k\) is well-conditioned.
Consequently, it is more reasonable to verify whether \(\vU_k\) is well-conditioned, i.e.,
\begin{equation*}
    \kappa(\vU_{k}) \leq \texttt{const},
\end{equation*}
where \(\texttt{const}\) denotes a constant that is independent of \(s\), \(m\), and \(\kappa(\bXX_k)\).
Notice that \(\Omega_{k} = \vU_{k}\trans \vU_{k}\) has already been computed in Line~\ref{line:bcgspipiro1s:sync} of each iteration.
Moreover, we can determine the dimension \(d\) to switch to \BCGSPIPIROtwo by verifying when the following condition holds:
\begin{equation} \label{eq:condition-switch}
    \texttt{const}^2 \lambda_{\min}(\Omega_d)
    = \texttt{const}^2 \sigmin^2(\vU_{d}) 
    \leq \norm{\vU_{d}}^2 
    = \lambda_{\max}(\Omega_d),
\end{equation}
where \(\lambda_{\min}\) and \(\lambda_{\max}\) denote the smallest and largest eigenvalue, respectively.
Observe that this operation only involves the small \(s\)-by-\(s\) matrix \(\Omega_{k}\), and hence it only requires an additional \(\bigO(s^3)\) operations per iteration and does not introduce additional synchronization points.

We summarize the adaptive approach combining \(\BCGSPIPIROtwo\) with \(\BCGSPIPIROone\) in Algorithm~\ref{alg:BCGSIROA+BCGSPIPIRO1s}.

\begin{algorithm}[htbp!]
	\caption{$[\bQQ, \RR] = \BCGSIROAPIPIROone(\bXX, \IOAnoarg, \IOonenoarg)$ \label{alg:BCGSIROA+BCGSPIPIRO1s}}
	\begin{algorithmic}[1]
            \State $[\bQQ_{d-1}, \RR_{d-1}] = \BCGSPIPIROone(\bXX, \IOAnoarg)$ employing~\eqref{eq:condition-switch} to find \(d\) such that \(\bXX_{d-1} = \bQQ_{d-1} \RR_{d-1}\)
		\For{$k = d, d+1, \ldots, p-1$}
		    \State $\SS_{1:k,k+1} =
                    \begin{bmatrix} \vZ_{k-1} \\
                    Y_{kk}\itrans \left( P_k - \YY_{1:k-1,k}\trans \vZ_{k-1} \right) \end{bmatrix}$ \label{line:1s2s:computeS}
                \State $[\vU_{k+1}, S_{k+1, k+1}] = \IOone{\vX_{k+1} - \bQQ_k \SS_{1:k,k+1}}$ \SComment{First synchronization} \label{line:1s2s:io1}
                \If{$k < p-1$}
                    \State $\begin{bmatrix} \YY_{1:k,k+1} & \vZ_{k} \\ \Omega_{k+1} & P_{k+1} \end{bmatrix} =
                    \begin{bmatrix}\bQQ_{k} & \vU_{k+1} \end{bmatrix}\trans \begin{bmatrix}\vU_{k+1} & \vX_{k+2} \end{bmatrix}$ \SComment{Second synchronization} \label{line:1s2s:sync}
                \Else
                    \State $\begin{bmatrix} \YY_{1:k,k+1}\\ \Omega_{k+1} \end{bmatrix} =
                    \begin{bmatrix}\bQQ_{k} & \vU_{k+1} \end{bmatrix}\trans \vU_{k+1}$
                \EndIf
                \State $Y_{k+1,k+1} = \chol(\Omega_{k+1} - \YY_{1:k,k+1}\trans \YY_{1:k,k+1})$
                \State $\vQ_{k+1} = (\vU_{k+1} - \bQQ_{k} \YY_{1:k,k+1}) Y_{k+1,k+1}^{-1}$ \label{line:1s2s:Qk}
                \State $\RR_{1:k,k+1} = \SS_{1:k,k+1} + \YY_{1:k,k+1} \diff{S_{k+1,k+1}}$
                \State $R_{k+1,k+1} = Y_{k+1,k+1} \diff{S_{k+1,k+1}}$
		\EndFor
		\State \Return{$\bQQ = [\vQ_1, \ldots, \vQ_p]$, $\RR = (R_{ij})$}
	\end{algorithmic}
\end{algorithm}

Based on Theorem~\ref{thm:pipiro2s}, \(\BCGSPIPIROone\) also satisfies a similar assumption as the \(\IOAnoarg\) required by Theorem~\ref{thm:pipiro2s}.
Thus, it is straightforward to generalize Theorem~\ref{thm:pipiro2s} to the following theorem for describing LOO of \(\BCGSIROAPIPIROone\).

\begin{theorem} \label{thm:bcgsiroapipiro1s}
    Let $\bQQhat$ and $\RRhat$ denote the computed results of Algorithm~\ref{alg:BCGSIROA+BCGSPIPIRO1s}.
    Assume that for all $\vX \in \mathcal{R}^{m \times s}$ with $\kappa(\vX) \leq \kappa(\bXX)$, the following hold for $[\vQhat, \Rhat] = \IOA{\vX}$ used in \BCGSPIPIROone:
    \begin{equation*}
        \begin{split}
            & \vX + \Delta \vX = \vQhat \Rhat,
            \quad \norm{\Delta \vX} \leq \bigO(\macheps) \norm{\vX}, \\
            & \norm{I - \vQhat\trans \vQhat} \leq \bigO(\macheps).
        \end{split}
    \end{equation*}
    Likewise, assume the following hold for $[\vQhat, \Rhat] = \IOone{\vX}$:
    \begin{equation*}
        \begin{split}
            & \vX + \Delta \vX = \vQhat \Rhat,
            \quad \norm{\Delta \vX} \leq \bigO(\macheps) \norm{\vX}, \\
            & \norm{I - \vQhat^T \vQhat} \leq \bigO(\macheps) \kappa^{\alpha_1}(\vX)
        \end{split}
    \end{equation*}
    with \(\alpha_1 \leq 1\).
    If $\bigO(\macheps) \kappa(\bXX) \leq 1/2$ is satisfied, then
    \begin{equation*}
        \bXX + \Delta\bXX = \bQQhat \RRhat,
        \quad \norm{\Delta\bXX} \leq \bigO(\macheps) \norm{\bXX},
    \end{equation*}
    and
    \begin{equation*}
        \norm{I - \bQQhat^T \bQQhat} \leq \bigO(\macheps).
    \end{equation*}
\end{theorem}

\subsection{\(s\)-step GMRES with low-sync BCGS}
In the $s$-step Arnoldi process, during the \(k\)-th iteration, we initially derive the Krylov basis \(B_{ks+1:(k+1)s}\) and subsequently execute the \(k\)-th iteration of a block orthogonalization method on \(\vX_{k+1} = M_L^{-1}AM_R^{-1}B_{ks+1:(k+1)s}\).
It should be noted that \(\vX_1 = \bmat{r& M_L^{-1}AM_R^{-1}B_{1:s}}\) consists of \(s+1\) columns for the Arnoldi process.
Notice that the \(k\)-th iteration of \BCGSPIPIROone and \BCGSPIPIROtwo also requires \(\vX_{k+2}\), which is generated by the last column of \(\vQ_{k+1}\) in the standard Arnoldi process.
However, to minimize synchronization points, the step involving \(\vX_{k+2}\) (refer to Line~\ref{line:bcgspipiro1s:sync} in Algorithm~\ref{alg:BCGSPIPIRO1S} or Line~\ref{line:bcgspipiro2s:sync} in Algorithm~\ref{alg:BCGSPIPIRO2S}) needs execution prior to forming \(\vQ_{k+1}\) (Line~\ref{line:Qk} in Algorithm~\ref{alg:BCGSPIPIRO1S} or Line~\ref{line:bcgspipiro2s:Qk} in Algorithm~\ref{alg:BCGSPIPIRO2S}), leaving \(\vQ_{k+1}\) unavailable for generating \(\vX_{k+2}\).
Hence, we instead use \(\vU_{k+1}\), which can be computed before Line~\ref{line:bcgspipiro1s:sync} of Algorithm~\ref{alg:BCGSPIPIRO1S} or Line~\ref{line:bcgspipiro2s:sync} of Algorithm~\ref{alg:BCGSPIPIRO2S}, to generate the basis, since \(\vU_{k+1} = \vQ_{k+1}\) in exact arithmetic.

In Algorithms~\ref{alg:sstep-Arnoldi} and~\ref{alg:sstep-GMRES}, we present the \(s\)-step Arnoldi process and \(s\)-step GMRES with \BCGSPIPIROtwo, respectively.
Employing \BCGSPIPIROone and \BCGSIROAPIPIROone in \(s\)-step GMRES follows a similar approach.

\begin{algorithm}[!tb]
\begin{algorithmic}[1]
    \caption{The \(s\)-step Arnoldi process with \BCGSPIPIROtwo \label{alg:sstep-Arnoldi}}
    \Require
    A matrix \(A \in \mathbb R^{m\times m}\), a right-hand side \(b \in \mathbb R^{m}\), an initial guess \(x^{(0)} \in \mathbb R^{m}\), a block size \(s\), a left-preconditioner \(M_L \in \mathbb R^{m\times m}\), and a right-preconditioner \(M_R \in \mathbb R^{m\times m}\).
    \Ensure
    A Krylov basis \(B_{1:ps}\) and a computed orthonormal basis \(Q_{1:ps+1}\) and upper triangular matrix \(R_{1:ps+1,1:ps+1}\) approximately satisfying \(\bmat{r& M_L^{-1} A M_R^{-1} B_{1:ps}} = Q_{1:ps+1} R_{1:ps+1,1:ps+1}\).

    \State \(r = M_L^{-1} (b - A x^{(0)})\) and \(R_{11} = \beta = \norm{r}\).
    \State \(Q_{1} = r/\beta\).
    \State \(B_{1:s} = \bmat{p_0(A) Q_{1} & p_1(A) Q_{1} & \dotsi & p_{s-1}(A) Q_{1}}\)
    \State \(X_{2:s+1} = M_L^{-1} A M_R^{-1} B_{1:s}\).
    \State $\SS = Q_1\trans X_{2:s+1}$
    \State $[U, S] = \IOone{X_{2:s+1} - Q_1 \SS}$
        \State \(B_{s+1:2s} = \bmat{p_0(A) U_{s} & p_1(A) U_{s} & \dotsi & p_{s-1}(A) U_{s}}\) \SComment{Generate the basis}
        \State \(X_{s+2:2s+1} = M_L^{-1} A M_R^{-1} B_{s+1:2s}\) \SComment{Apply preconditioners to the basis}
        \State {$\begin{bmatrix} \YY & Z \\ \Omega & P \end{bmatrix} = \begin{bmatrix}Q_{1} & U \end{bmatrix}\trans \begin{bmatrix} U & X_{s+2:2s+1} \end{bmatrix}$}
        \State $Y = \chol(\Omega - \YY\trans \YY)$
        \State $Q_{2:s+1} = (U - Q_{1} \YY) Y^{-1}$
        \State $R_{1,2:s+1} = \SS + \YY S$
        \State $R_{2:s+1,2:s+1} = Y S$
        \For{$k = 2, \ldots, p$}
            \State $\SS =
                \begin{bmatrix} Z \\
                Y\itrans \left( P - \YY\trans Z \right) \end{bmatrix}$
            \State $[U, S] = \IOone{X_{ks+2:(k+1)s+1} - Q_{1:ks+1} \SS}$ 
            \State \(B_{ks+1:(k+1)s} = \bmat{p_0(A) U_{s} & p_1(A) U_{s} & \dotsi & p_{s-1}(A) U_{s}}\) \SComment{Generate the basis}
            \State \(X_{ks+2:(k+1)s+1} = M_L^{-1} A M_R^{-1} B_{ks+1:(k+1)s}\) \SComment{Apply preconditioners to the basis}
            \State $\begin{bmatrix} \YY & Z \\ \Omega & P \end{bmatrix} =
                \begin{bmatrix} Q_{1:ks+1} & U\end{bmatrix}\trans \begin{bmatrix} U & X_{ks+2:(k+1)s+1} \end{bmatrix}$
            \State $Y = \chol(\Omega - \YY\trans \YY)$
            \State $Q_{(k-1)s+2:ks+1} = (U - Q_{1:(k-1)s+1} \YY) Y^{-1}$
            \State $\RR_{1:(k-1)s+1,(k-1)s+2:ks+1} = \SS + \YY \diff{S}$
            \State $R_{(k-1)s+2:ks+1, (k-1)s+2:ks+1} = Y \diff{S}$
        \EndFor
\end{algorithmic}
\end{algorithm}

\begin{algorithm}[!tb]
\begin{algorithmic}[1]
    \caption{The \(s\)-step GMRES algorithm \label{alg:sstep-GMRES}}
    \Require
    A matrix \(A \in \mathbb R^{m\times m}\), a right-hand side \(b \in \mathbb R^{m}\), an initial guess \(x^{(0)} \in \mathbb R^{m}\), a block size \(s\), a left-preconditioner \(M_L \in \mathbb R^{m\times m}\), and a right-preconditioner \(M_R \in \mathbb R^{m\times m}\).
    \Ensure
    A computed solution \(x \in \mathbb R^{m}\) approximating the solution of \(A x = b\).

    \State \(r = M_L^{-1} (b - A x^{(0)})\) and \(\beta = \norm{r}\).
    \State \(Q_{1} = r/\beta\).
    \For{\(k = 1:m/s\)}
        \State Perform the \(k\)-th step of the \(s\)-step Arnoldi process (e.g., Algorithm~\ref{alg:sstep-Arnoldi}) to obtain the basis \(B_{1:ks}\), the preconditioned basis \(Z_{1:ks} = M_R^{-1} B_{1:ks}\), the orthonormal matrix \(Q_{1:ks+1}\), and the upper triangular matrix \(R_{1:ks+1}\) satisfying \(\bmat{r & W_{1:ks}} = Q_{1:ks+1} R_{1:ks+1, 1:ks+1}\) with \(W_{1:ks} = M_L^{-1} A Z_{1:ks}\).
        \State \(H_{1:(k-1)s+1, 1:(k-1)s} = R_{1:(k-1)s+1, 2:(k-1)s+1}\).
        \State Compute the QR factorization \(H_{1:ks+1, (k-1)s+1:ks} = G_{1:ks+1} T_{1:ks+1, (k-1)s+1:ks}\) by Givens rotations with an orthogonal matrix \(G_{1:ks+1}\), based on \(H_{1:(k-1)s+1, 1:(k-1)s} = G_{1:(k-1)s+1} T_{1:(k-1)s+1, 1:(k-1)s}\).
        \If{the stopping criterion is satisfied}
             \State Solve the triangular system \(T_{1:ks} y^{(k)} = \beta G_{1, 1:ks}\trans\) to obtain \(y^{(k)} \in \mathbb R^{ks}\).
             \State \Return \(x = x^{(k)} = x^{(0)} + Z_{1:ks} y^{(k)}\).
        \EndIf
    \EndFor
\end{algorithmic}
\end{algorithm}

The \(s\)-step GMRES algorithm with \BCGSPIPIROtwo or \BCGSIROAPIPIROone can also be proven to be backward stable just as was done for the \(s\)-step GMRES algorithm with \BCGSIRO~\cite{CM2024}.
For presenting the backward error analysis, we first prove the following corollary about \BCGSPIPIROtwo and \BCGSIROAPIPIROone.

\begin{corollary} \label{cor:orth}
    Let $\bQQhat$ and $\RRhat$ denote the computed results of Algorithm~\ref{alg:BCGSPIPIRO2S} or~\ref{alg:BCGSIROA+BCGSPIPIRO1s}.
    Assume that for all $\vX \in \mathcal{R}^{m \times s}$ with $\kappa(\vX) \leq \kappa(\bXX)$, the following hold for $[\vQhat, \Rhat] = \IOA{\vX}$:
    \begin{equation*}
        \begin{split}
            & \vX + \Delta \vX = \vQhat \Rhat,
            \quad \norm{\Delta \vX} \leq \bigO(\macheps) \norm{\vX}, \\
            & \norm{I - \vQhat\trans \vQhat} \leq \bigO(\macheps).
        \end{split}
    \end{equation*}
    Likewise, assume the following hold for $[\vQhat, \Rhat] = \IOone{\vX}$:
    \begin{equation} \label{eq:cor:IO1assump}
        \begin{split}
            & \vX + \Delta \vX = \vQhat \Rhat,
            \quad \norm{\Delta \vX} \leq \bigO(\macheps) \norm{\vX}, \\
            & \norm{I - \vQhat^T \vQhat} \leq \bigO(\macheps) \kappa^{\alpha_1}(\vX)
        \end{split}
    \end{equation}
    with \(\alpha_1 \leq 1\).
    If $\bigO(\macheps) \kappa(\bXX) \leq 1/2$ is satisfied, then there exists an orthonormal matrix \(\bQQtil\) such that
    \begin{equation} \label{eq:cor:res}
        \bXX + \Delta\widetilde\bXX = \bQQtil \RRhat,
        \qquad \norm{\Delta\widetilde X_i} \leq \bigO(\macheps) \norm{X_i}, \quad \forall i = 1, 2, \dotsc, n
    \end{equation}
    and
    \begin{equation} \label{eq:cor:loo}
        \norm{I - \bQQhat^T \bQQhat} \leq \bigO(\macheps).
    \end{equation}
\end{corollary}

\medskip
\begin{proof}
    Equation~\eqref{eq:cor:loo} can be directly derived from Theorems~\ref{thm:pipiro2s} and~\ref{thm:bcgsiroapipiro1s}, which implies that
    \begin{equation} \label{eq:cor-proof:sigminbQQhat}
        \sigmin(\bQQhat) \geq 1 - \bigO(\macheps).
    \end{equation}
    Then we only need to show~\eqref{eq:cor:res}.
    Notice that every operation in \BCGSPIPIROtwo and \BCGSIROAPIPIROone is column-wise backward stable for \(\bXX\).
    Together with Theorems~\ref{thm:pipiro2s} and~\ref{thm:bcgsiroapipiro1s}, we have
    \begin{equation} \label{eq:cor-proof:res}
        \bXX + \Delta\bXX = \bQQhat \RRhat,
        \qquad \norm{\Delta X_i} \leq \bigO(\macheps) \norm{X_i}, \quad \forall i = 1, 2, \dotsc, n.
    \end{equation}
    Combined~\eqref{eq:cor-proof:sigminbQQhat} with~\eqref{eq:cor-proof:res}, it follows that
    \begin{equation*}
        \sigmin(\bQQhat) \norm{\Rhat_i} \leq \norm{\bQQhat \Rhat_i} \leq \norm{X_i} + \norm{\Delta X_i} \leq (1 + \bigO(\macheps)) \norm{X_i},
    \end{equation*}
    which means that
    \begin{equation} \label{eq:cor-proof:normRhati}
        \norm{\Rhat_i} \leq \frac{(1 + \bigO(\macheps)) \norm{X_i}}{\sigmin(\bQQhat)} \leq \frac{1 + \bigO(\macheps)}{1 - \bigO(\macheps)} \norm{X_i}.
    \end{equation}
    Considering the polar decomposition of \(\bQQhat\), i.e., \(\bQQhat = \bQQtil H\), where \(\bQQtil\) is orthonormal and \(H\) is symmetric positive semidefinite, by~\cite[Lemma~5.1]{Higham1994} and~\eqref{eq:cor:loo}, we obtain
    \[
        \norm{\bQQhat - \bQQtil} \leq \norm{I - \bQQhat^T \bQQhat} \leq \bigO(\macheps).
    \]
    Furthermore, together with~\eqref{eq:cor-proof:res} and~\eqref{eq:cor-proof:normRhati}, we conclude the proof of~\eqref{eq:cor:res} because
    \begin{equation}
        \bQQtil \RRhat = \bQQhat \RRhat + (\bQQtil - \bQQhat) \RRhat = \bXX + \underbrace{\Delta\bXX + (\bQQtil - \bQQhat) \RRhat}_{=: \Delta\widetilde\bXX}
    \end{equation}
    with \(\norm{\Delta\widetilde X_i} \leq \norm{\Delta X_i} + \norm{\bQQtil - \bQQhat} \norm{\Rhat_i} \leq \bigO(\macheps) \norm{X_i}\) for any \(i = 1, 2, \dotsc, n\).
\end{proof}

We now present the backward error analysis of the \(s\)-step GMRES algorithm with \BCGSPIPIROtwo or \BCGSIROAPIPIROone in the following lemma derived directly by using Corollary~\ref{cor:orth} and~\cite[Theorem 1]{CM2024}.

\begin{lemma} \label{lem:sstep-GMERS}
    Assume that the Krylov basis \(\widehat{B}_{1:ks}\) and the approximate solution \(\hat{x}\) are computed by the \(s\)-step GMRES algorithm with \BCGSPIPIROtwo or \BCGSIROAPIPIROone orthogonalization, in order to solve the linear system \(Ax = b\).
    Also, assume that the computed results \(\widehat{Z}_{1:ks}\) and \(\widehat{W}_{1:ks}\) satisfy
    \begin{align}
        \widehat{Z}_{1:ks} & = M_R^{-1} \hat B_{1:ks} + \Delta Z_{1:ks}, \quad \norm{\Delta Z_{j}} \leq \bigO(\macheps) \norm{M_R^{-1}} \norm{\hat{B}_{j}}, \label{eq:sstep-GMRES-barB}\\
        \widehat W_{1:ks} & = M_L^{-1} A \widehat{Z}_{1:ks} + \Delta D_{1:ks}, \quad \norm{\Delta D_{j}} \leq \bigO(\macheps) \norm{M_L^{-1}} \normF{A} \norm{\widehat{Z}_{j}}
    \end{align}
    for any \(j \leq ks\).
    There exists \(\kn = \ki s + \kj\) such that
    \begin{align}
        & \sigmin\Big(\bmat{\phi \widehat{r}& \widehat{W}_{1:\kn} D_{1:\kn}^{-1}}\Big) \leq \bigO(\macheps) \normF{\bmat{\phi \widehat{r}& \widehat{W}_{1:\kn} D_{1:\kn}^{-1}}} \label{eq:lem:LS:assump-rW}
    \end{align}
    holds for any \(\phi \geq 0\).
    If it also holds that
    \begin{equation*}
        \frac{\bigO(\macheps) \kappa(M_L) \kappa(A) \kappa(\tilde Z_{1:\kn})}{1 - \bigO(\macheps) \kappa(M_L) \kappa(A) \kappa(\tilde Z_{1:\kn})} \leq 1,
    \end{equation*}
    then
    \begin{equation*}
        \begin{split}
            \frac{\norm{b - A\widehat{x}^{(\ki)}}}{\norm{b} + \normF{A} \norm{\widehat x^{(\ki)}}}
            \leq{}& \frac{\bigO(\macheps) \kappa(M_L) \kappa (\tilde B_{1:\kn})}{1 - \bigO(\macheps) \kappa (\tilde B_{1:\kn})},
        \end{split}
    \end{equation*}    
    where \(\widehat B_{1:\kn} = \tilde B_{1:\kn} D_{1:\kn}\) and \(\widehat Z_{1:\kn} = \tilde Z_{1:\kn} D_{1:\kn}\) with a positive definite diagonal matrix \(D_{1:\kn}\).
\end{lemma}
\section{Numerical Experiments}
\label{sec:experiments}
In this section, we first compare the relative residual and LOO of \(\BCGSIRO\), \(\BCGSPIPIROone\), \(\BCGSPIPIROtwo\), \(\BCGSIROAPIPIROone\), \(\BCGSPIPIRO\), and \(\BCGSIROAone\), where \BCGSIROAone is the previous one-sync reorthogonalized BCGS method equivalent to \texttt{BCGS+LS} proposed in~\cite{YTHBSE2020}.
Subsequently, we evaluate the impact of employing these BCGS variants within \(s\)-step GMRES.

In~\eqref{eq:condition-switch}, we set ``\(\texttt{const}\)'' as \(\sqrt{3}\) for switching to \BCGSPIPIROtwo in \BCGSIROAPIPIROone, utilizing \HouseQR for both \(\IOAnoarg\) and \(\IOonenoarg\).
Similarly, for \BCGSIRO, \HouseQR is employed as \(\IOonenoarg\) and \(\IOtwonoarg\).
This gives similar numerical results as we would expect using TSQR.
For \(s\)-step GMRES, the stopping criterion is defined by
\[
\norm{b - A x^{(k)}} \leq 10^{-12} \bigl(\normF{A} \norm{x^{(k)}} + \norm{b}\bigr),
\]
where \(x^{(k)}\) represents the approximate solution at the \(k\)-th iteration.
All experiments are conducted using MATLAB R2024a.

\subsection{Comparison between different variants of BCGS}
We employ the MATLAB code suite \texttt{BlockStab}%
\footnote{https://github.com/katlund/BlockStab}
to compare different variants of BCGS across four matrix classes provided in \texttt{BlockStab}, namely \default, \glued, \monomial, \piled, as detailed and evaluated in~\cite{CLMO2024-ls,CLMO2024-P}.
These matrices are designed as commonly used and challenging test cases for BCGS variants.
\emph{In the figures, we display only the LOO plots, as the backward errors \(\norm{\bXX-\hat{\bQQ}\hat{\RR}}/\norm{\bXX}\) for each method are \(\bigO(\macheps)\).}

\begin{figure}[!tb]
\includegraphics[width=0.8\textwidth]{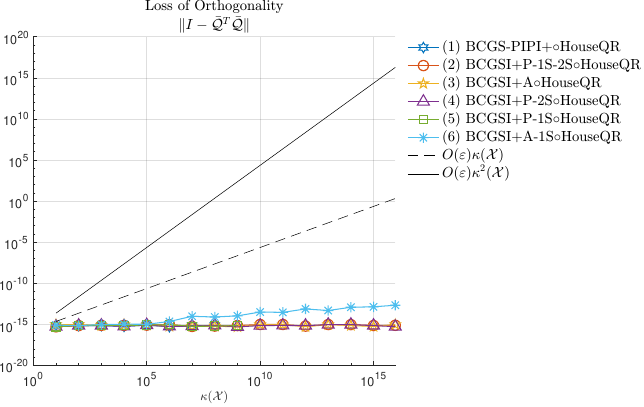}
\centering
\caption{Comparison between \BCGSPIPIROone, \BCGSPIPIROtwo, \BCGSIROAone, \BCGSPIPIRO,  \BCGSIROAPIPIROone, and \BCGSIRO on a class of \default matrices.}
\label{fig:default}
\end{figure}

\begin{figure}[!tb]
\includegraphics[width=0.8\textwidth]{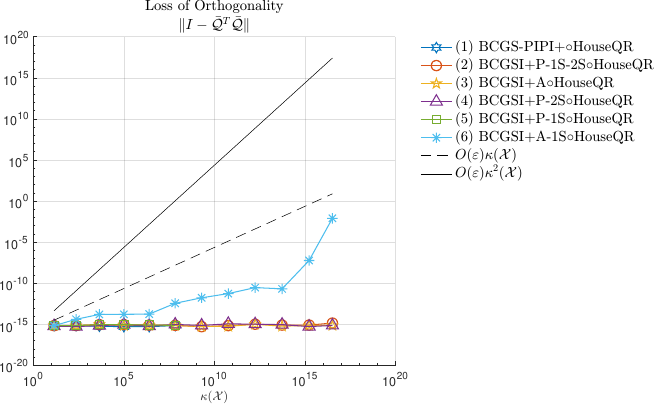}
\centering
\caption{Comparison between \BCGSPIPIROone, \BCGSPIPIROtwo, \BCGSIROAone, \BCGSPIPIRO,  \BCGSIROAPIPIROone, and \BCGSIRO on a class of \glued matrices.}
\label{fig:glued}
\end{figure}

\begin{figure}[!tb]
\includegraphics[width=0.8\textwidth]{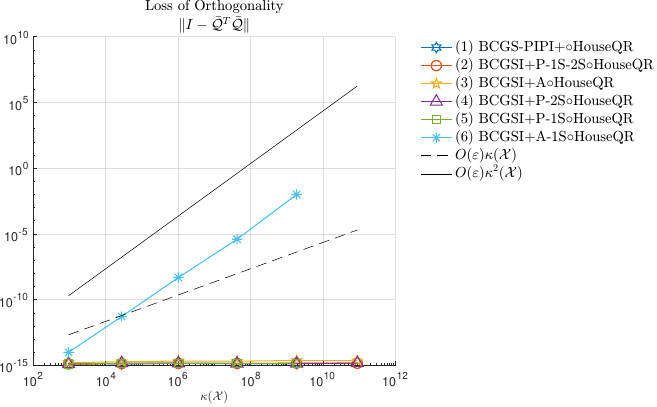}
\centering
\caption{Comparison between \BCGSPIPIROone, \BCGSPIPIROtwo, \BCGSIROAone, \BCGSPIPIRO,  \BCGSIROAPIPIROone, and \BCGSIRO on a class of \monomial matrices.}
\label{fig:monomial}
\end{figure}

\begin{figure}[!tb]
\includegraphics[width=0.8\textwidth]{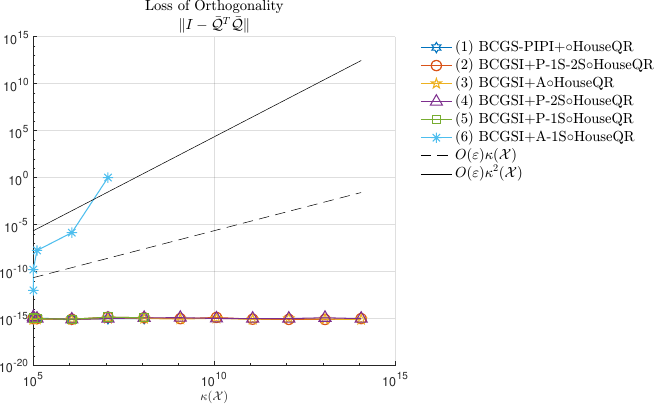}
\centering
\caption{Comparison between \BCGSPIPIROone, \BCGSPIPIROtwo, \BCGSIROAone, \BCGSPIPIRO,  \BCGSIROAPIPIROone, and \BCGSIRO on a class of \piled matrices.}
\label{fig:piled}
\end{figure}

Figures~\ref{fig:default}--\ref{fig:piled} show the behavior of \(\BCGSPIPIROone\), \(\BCGSPIPIROtwo\), \(\BCGSIROAPIPIROone\), \(\BCGSPIPIRO\), \(\BCGSIROAone\), and \(\BCGSIRO\) (i.e., \BCGSIROA with \(\IOAnoarg=\HouseQR\)) over various matrix classes.
It is observed that \(\BCGSIROAone\) is unstable, particularly for \monomial and \piled matrices.
In contrast, our one-sync variant \(\BCGSPIPIROone\) is as stable as \(\BCGSPIPIRO\) and exhibits \(\bigO(\macheps)\) LOO for all test matrices with \(\bigO(\macheps) \kappa^2(\bXX) \leq 1\).
Furthermore, our two-sync variant \(\BCGSPIPIROtwo\) and \(\BCGSIROAPIPIROone\) can be used with all matrices satisfying \(\bigO(\macheps) \kappa(\bXX) \leq 1\), maintaining the same LOO as \(\BCGSIRO\) while reducing the number of synchronization points during the iterations.

\subsection{\(s\)-step GMRES with different variants of BCGS}
Here we evaluate the behavior of \(s\)-step GMRES with several BCGS variants, specifically \BCGSPIPIROone, \BCGSPIPIROtwo, \BCGSIROAPIPIROone, and \BCGSIRO, for solving the linear system \(Ax = b\).
To construct the linear systems, we use sparse matrices shown in Table~\ref{table:testmat} from the SuiteSparse Matrix Collection%
\footnote{\url{https://sparse.tamu.edu}}
and the Matrix Market%
\footnote{\url{https://math.nist.gov/MatrixMarket}}
as \(A\).
The right-hand side vector \(b\) consists entirely of ones, and the initial guess \(x^{(0)}\) is set to zero.

\begin{figure}[!tb]
\centering
\includegraphics[width=\textwidth]{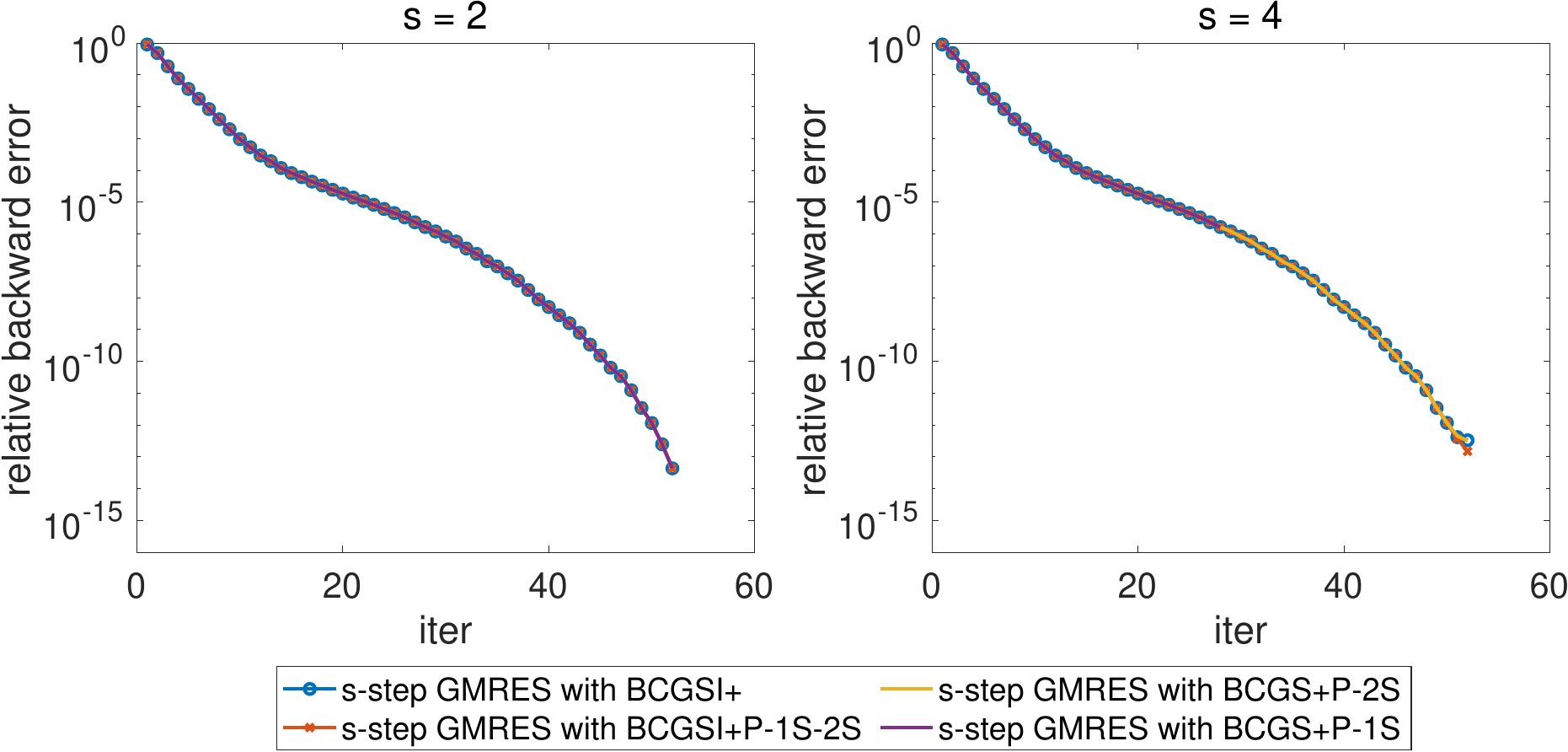}
\caption{Comparison between \(s\)-step GMRES using \BCGSPIPIROone, \BCGSPIPIROtwo, \BCGSIROAPIPIROone, and \BCGSIRO as orthogonalization on \texttt{fs\_760\_1} matrix.}
\label{fig:fs760-1}
\end{figure}

Figure~\ref{fig:fs760-1} illustrates that the behavior of \(s\)-step GMRES of a well-conditioned matrix is comparable across \BCGSPIPIROtwo, \BCGSIROAPIPIROone, and \BCGSIRO.
In the figures and the table, the relative backward error means that \(\norm{b - A \hat x^{(k)}}/ \bigl(\normF{A} \norm{\hat x^{(k)}} + \norm{b}\bigr)\).
Note that for $s=4$, the $s$-step GMRES method using \BCGSPIPIROtwo fails after reaching the level $O(\sqrt{\macheps})$, as predicted by the theory.
As shown in Table~\ref{table:info-fs760-1}, \(s\)-step GMRES using \BCGSPIPIROtwo and \BCGSIROAPIPIROone for orthogonalization can also achieve the same level of backward error as using \BCGSIRO, but requires fewer synchronization points.

\begin{figure}[!tb]
\centering
\includegraphics[width=\textwidth]{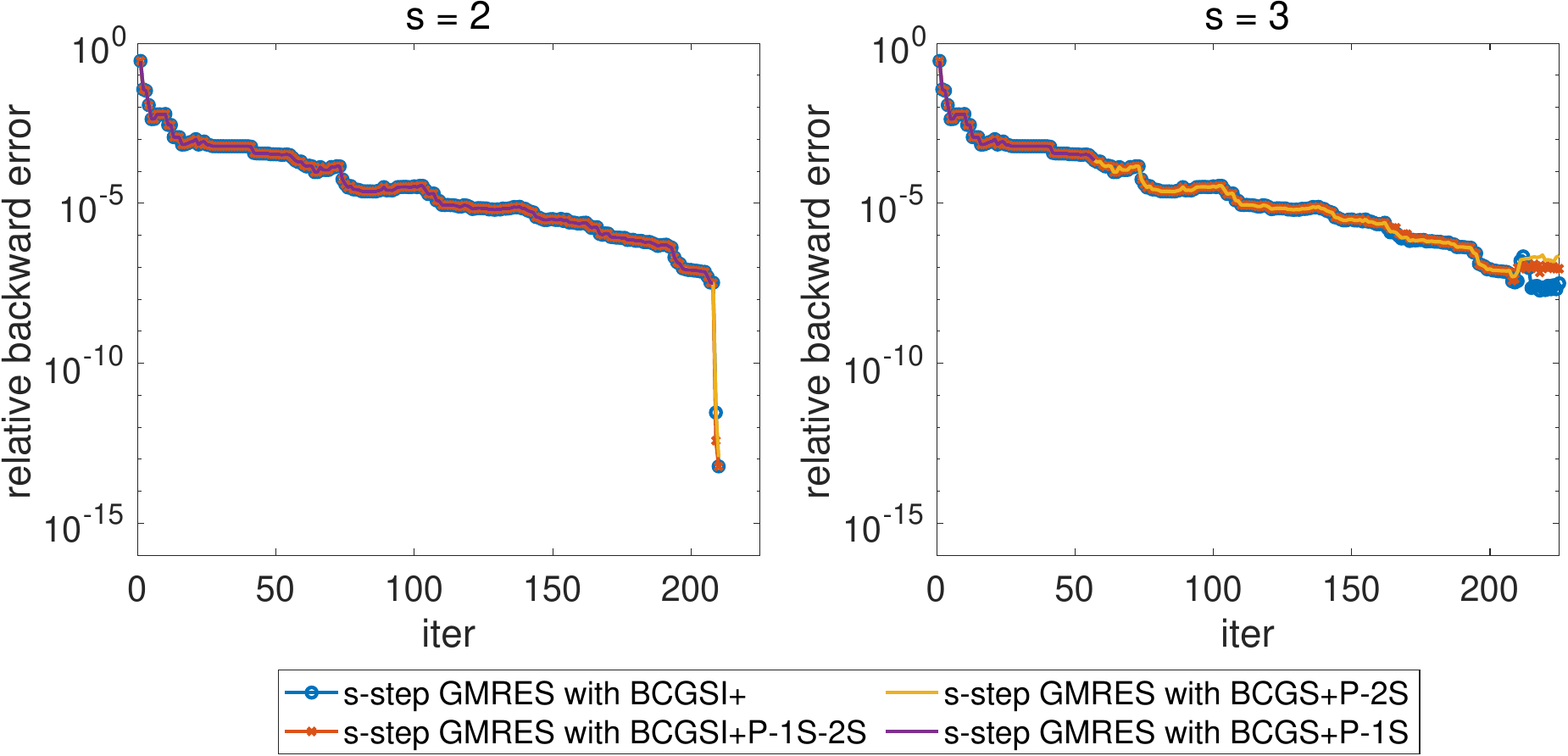}
\caption{Comparison between \(s\)-step GMRES using \BCGSPIPIROone, \BCGSPIPIROtwo, \BCGSIROAPIPIROone, and \BCGSIRO as orthogonalization on \texttt{impcol\_e} matrix.}
\label{fig:impcole}
\end{figure}

\begin{figure}[!tb]
\centering
\includegraphics[width=\textwidth]{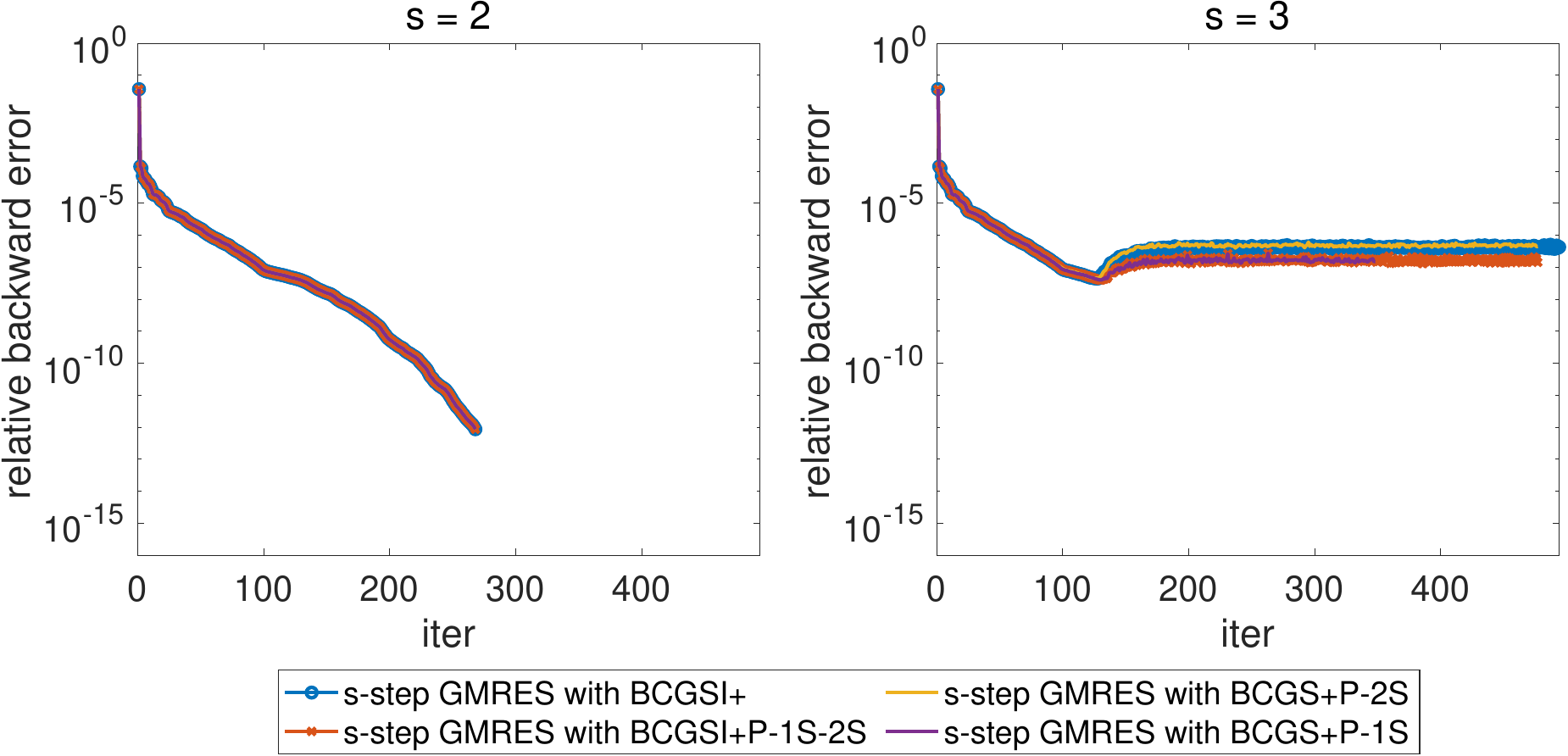}
\caption{Comparison between \(s\)-step GMRES using \BCGSPIPIROone, \BCGSPIPIROtwo, \BCGSIROAPIPIROone, and \BCGSIRO as orthogonalization on \texttt{494bus} matrix.}
\label{fig:494bus}
\end{figure}

\begin{figure}[!tb]
\centering
\includegraphics[width=\textwidth]{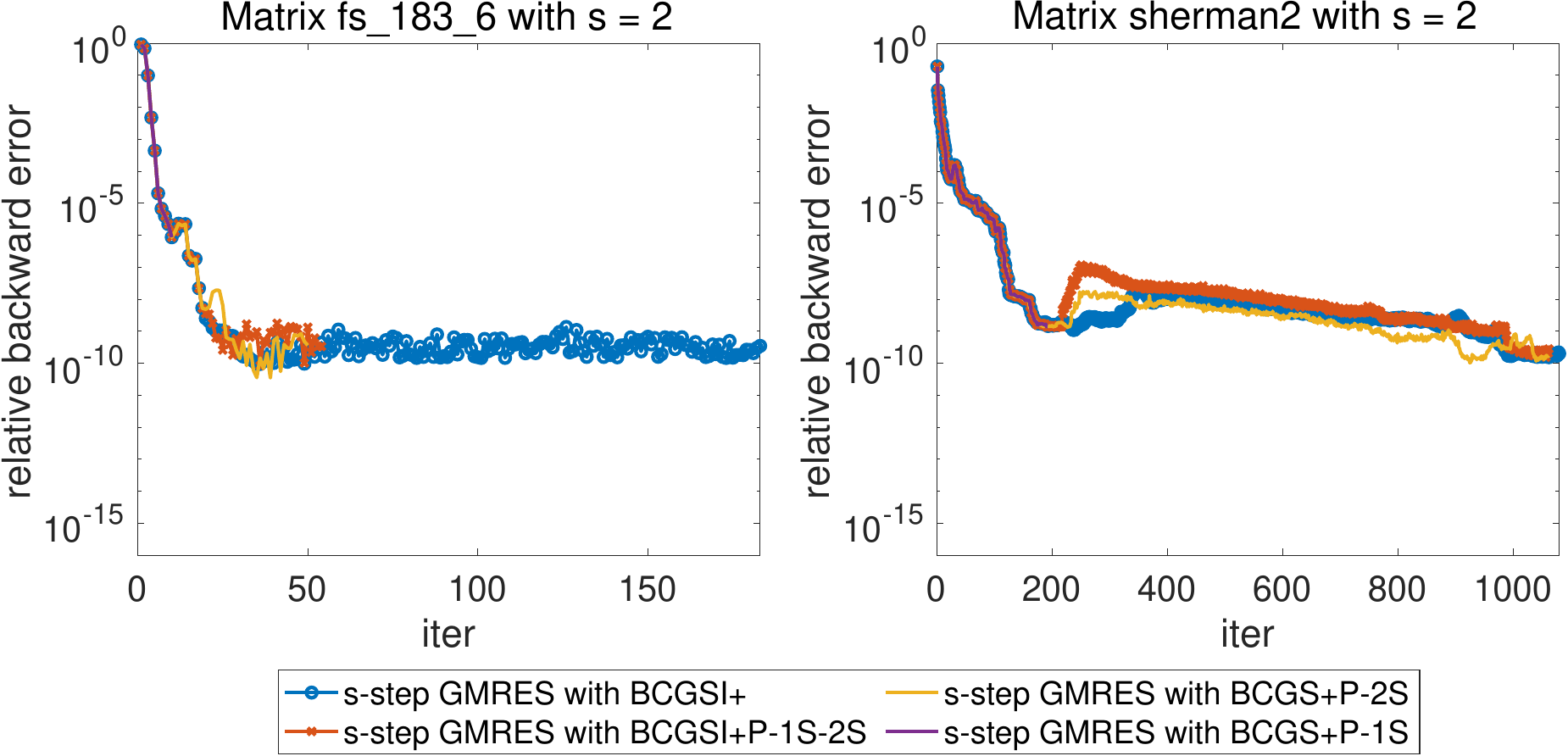}
\caption{Comparison between \(s\)-step GMRES using \BCGSPIPIROone, \BCGSPIPIROtwo, \BCGSIROAPIPIROone, and \BCGSIRO as orthogonalization on \texttt{sherman2} and \texttt{fs\_183\_6} matrices.}
\label{fig:sherman2-fs1836}
\end{figure}

Additionally, we test more matrices with larger condition numbers as shown in Figures~\ref{fig:impcole}--\ref{fig:sherman2-fs1836}.
For each of these cases, \(s\) is chosen to ensure \(\bigO(\macheps) \kappa(B_{(k-1)s+1:ks}) \leq 1\) in the \(k\)-th iteration.
The figures demonstrate that the \(s\)-step GMRES method employing \BCGSPIPIROtwo, \BCGSIROAPIPIROone, and \BCGSIRO behaves similarly in most scenarios.
In particular, for the matrix \texttt{fs\_183\_6}, \(s\)-step GMRES with \BCGSPIPIROtwo and \BCGSIROAPIPIROone halts after approximately \(54\) iterations because the condition number of \(\bXX = \bmat{r& W}\) exceeds \(1/\bigO(\macheps)\), which causes the Cholesky factorization to fail.
Fortunately, this situation indicates that \(s\)-step GMRES has already reached the key dimension \(\kn\) as defined by~\eqref{eq:lem:LS:assump-rW} in Lemma~\ref{lem:sstep-GMERS}.
Continuing the iterations beyond this point does not improve the backward error. 
Consequently, these figures indicate that our low-synchronization BCGS algorithms \BCGSPIPIROtwo and \BCGSIROAPIPIROone are feasible for \(s\)-step GMRES without compromising the convergence rate or backward error.

It is important to acknowledge that \(s\)-step GMRES with our low-synchronization BCGS and \BCGSIRO will behave differently when the condition \(\bigO(\macheps) \kappa(B_{(k-1)s+1:ks}) \leq 1\) is not met.
In such cases, while \BCGSIRO can still generate an orthonormal basis, the low-synchronization BCGS algorithms, which employ Cholesky-based QR as the second intraorthogonalization, will fail.
Nevertheless, choosing a larger \(s\) with \(\bigO(\macheps) \kappa(B_{(k-1)s+1:ks}) > 1\) is often pointless because certain columns of \(B_{(k-1)s+1:ks}\) are ineffective for convergence if \(B_{(k-1)s+1:ks}\) is numerically singular. Additionally, in practice, it is easy to switch to employ \HouseQR or \TSQR in Line~\ref{line:bcgspipiro2s:Qk} of Algorithm~\ref{alg:BCGSPIPIRO2S} and Line~\ref{line:1s2s:Qk} of Algorithm~\ref{alg:BCGSIROA+BCGSPIPIRO1s} at the cost of adding an additional synchronization point when it is hard to choose a proper \(s\) in advance.

\begin{table}[htb!]
\centering
\caption{Properties of test matrices: the condition number in this table is estimated by the MATLAB command \texttt{svd}.}
\label{table:testmat}
\begin{tabular}{cccc}
\hline
Name& Size& Condition number& Source \\
\hline
\texttt{fs\_760\_1}      & \hphantom{0}760 & \(5.49\times 10^3\)\hphantom{0} &SuiteSparse Matrix Collection \\
\texttt{494bus}      & \hphantom{0}494 & \(2.42\times 10^6\)\hphantom{0} &SuiteSparse Matrix Collection \\
\texttt{impcol\_e}      & \hphantom{0}494 & \(7.10\times 10^6\)\hphantom{0} &SuiteSparse Matrix Collection \\
\texttt{fs1836}        & \hphantom{0}183 & \(1.74\times 10^{11}\) &Matrix Market \\
\texttt{sherman2}      & 1,080 & \(9.64\times 10^{11}\) &Matrix Market \\
\hline
\end{tabular}
\end{table}

\begin{table}[htb!]
\centering
\caption{Comparison of \(s\)-step GMRES incorporating \BCGSPIPIROone, \BCGSPIPIROtwo, \BCGSIROAPIPIROone, and \BCGSIRO for orthogonalization applied to the \texttt{fs\_760\_1} matrix. The term ``sync-points'' represents the total number of synchronization points during the Arnoldi process, with the exclusion of synchronization points used for calculating the first column of the orthonormal basis.
For \(s\)-step GMRES with \BCGSIROAPIPIROone, we also present in parentheses the iterations performed by \BCGSPIPIROone and \BCGSPIPIROtwo, respectively.
This implies that \(d-1\) in Algorithm~\ref{alg:BCGSIROA+BCGSPIPIRO1s} is \(26 (=52/s)\) and \(7 (=28/s)\), respectively, for \(s=2\) and \(s=4\).
Note that \(s\)-step GMRES with \BCGSPIPIROone with \(s=4\) does not converge, since a quantity becomes  \texttt{NaN} after 28 iterations.}
\vspace{2mm}
\label{table:info-fs760-1}
\begin{tabular}{ccccc}
\hline
Orthogonalization     & \(s\)& Relative backward error & Iterations & Sync-points \\
\hline
\BCGSIRO              & \(2\) & \(4.36\times 10^{-14}\) & 52\hphantom{ (24+28)} & 104\\
\BCGSIROAPIPIROone    & \(2\) & \(4.36\times 10^{-14}\) & 52 (52+0)\hphantom{0} & 26\hphantom{0}\\ 
\BCGSPIPIROtwo        & \(2\) & \(4.36\times 10^{-14}\) & 52\hphantom{ (24+28)} & 52\hphantom{0}\\
\BCGSPIPIROone        & \(2\) & \(4.36\times 10^{-14}\) & 52\hphantom{ (24+28)} & 26\hphantom{0}\\ \hline
\BCGSIRO              & \(4\) & \(5.75\times 10^{-13}\) & 52\hphantom{ (24+28)} & 52\hphantom{0}\\
\BCGSIROAPIPIROone    & \(4\) & \(1.69\times 10^{-13}\) & 52 (28+24) & 20\hphantom{0}\\
\BCGSPIPIROtwo        & \(4\) & \(2.21\times 10^{-13}\) & 52\hphantom{ (24+28)} & 26\hphantom{0}\\
\BCGSPIPIROone        & \(4\) & \(1.66\times 10^{-6\hphantom{0}}\)  & 28\hphantom{ (24+28)} & 7\hphantom{00}\\
\hline
\end{tabular}
\end{table}

\section{Conclusions}
\label{sec:conclusions}
For many years, researchers have sought a stable BCGS variant that requires only one synchronization point per iteration~\cite{CLMO2024-ls, SLAYT2021, YTHBSE2020}.
Existing BCGS algorithms, such as \BCGSIROAone in \cite{CLMO2024-ls} or \texttt{BCGS+LS} in~\cite{YTHBSE2020}, require just one synchronization per iteration but lack stability, whereas algorithms like \BCGSPIPIRO are stable but necessitate two synchronization points per iteration.
In this work, we propose a new one-sync reorthogonalized BCGS method, \BCGSPIPIROone, and show that it can achieve \(\bigO(\macheps)\) loss of orthogonality under the condition \(\bigO(\macheps) \kappa^2(\bXX) \leq 1/2\).
Our analysis of \BCGSPIPIROone reveals that this constraint arises from employing Cholesky-based QR for the first intraorthogonalization.
Therefore, we further propose a two-sync reorthogonalized BCGS, called \BCGSPIPIROtwo, which maintains \(\bigO(\macheps)\) loss of orthogonality under the improved condition \(\bigO(\macheps) \kappa(\bXX) \leq 1/2\), by incorporating a stable method as the first intraorthogonalization such as \HouseQR and \TSQR at the cost of adding an additional synchronization point per iteration.

Next, we discuss a crucial application of these low-synchronization BCGS variants, i.e., \(s\)-step GMRES, which requires that the orthogonalization method has \(\bigO(\macheps)\) loss of orthogonality under the condition \(\bigO(\macheps) \kappa^\alpha(\bXX) \leq 1\) with \(\alpha = 0\) or \(1\) to achieve a backward error at the level \(\bigO(\macheps)\).
This means that \BCGSPIPIROone is unsuitable for orthogonalization in \(s\)-step GMRES.
Thus, we propose an adaptive technique to combine \BCGSPIPIROone and \BCGSPIPIROtwo.
This approach also maintains \(\bigO(\macheps)\) loss of orthogonality under the condition \(\bigO(\macheps) \kappa(\bXX) \leq 1/2\), while requiring as few synchronization points as possible.
Furthermore, we analyze the backward stability of \(s\)-step GMRES utilizing these low-synchronization variants, namely \BCGSPIPIROtwo and \BCGSIROAPIPIROone.
Both methods can be proven to be as stable as \(s\)-step GMRES with \BCGSIRO.
Numerical experiments verify our theoretical results by comparing the loss of orthogonality of different BCGS variants and their application in \(s\)-step GMRES.
\section*{Acknowledgments}
Both authors are supported by the European Union (ERC, inEXASCALE, 101075632). Views and opinions expressed are those of the authors only and do not necessarily reflect those of the European Union or the European Research Council. Neither the European Union nor the granting authority can be held responsible for them. Both authors additionally acknowledge support from the Charles University Research Centre program No. UNCE/24/SCI/005.

\bibliographystyle{abbrvurl}
\bibliography{mybib}

\end{document}